\documentclass[11pt]{article}

\usepackage[latin1]{inputenc}  
\usepackage{amssymb,amsmath,latexsym,enumerate,multirow,amsthm}
\usepackage{graphicx,psfrag} 
\usepackage[tight]{subfigure} 
\usepackage{makeidx}
\usepackage{hyperref}
\usepackage{color}
\usepackage[english]{babel}

\newtheorem{algorithm}{Algorithm}
\newtheorem{testproblem}{Test Problem}
\newtheorem*{keywords}{Keyword}

\DeclareMathOperator{\dist}{dist}
\DeclareMathOperator{\avg}{avg}

\def\figpath{Figure_eps/}

\def\RR{\mathbb{R}}

\newcommand{\NN}{\mathbb{N}}

\begin{document}

\title{Improved Stencil Selection for Meshless Finite Difference Methods  in 3D}

\author{
Oleg Davydov\thanks{Department of Mathematics, University of Giessen, Arndtstrasse 2, 35392 Giessen, Germany,
{\tt oleg.davydov@math.uni-giessen.de}.},\quad
Dang Thi Oanh\thanks{%
Information Technology Department,
Ministry of Education and training,
35 Dai Co Viet, Hai Ba Trung District, Hanoi, Vietnam, {\tt dtoanh@moet.gov.vn; dtoanh@ictu.edu.vn}.
The work of this author was supported 
by a Natural Science Research Project of the Ministry of Education and
Training under grant number B2019-TNA-03.T.}\quad 
and
 Ngo Manh Tuong\thanks{ Department of Basic Sciences,
Thai Nguyen University of Information $\&$ Communication Technology,
Quyet Thang Ward, Thai Nguyen City, Vietnam, {\tt nmtuong@ictu.edu.vn}.}
 }
\maketitle

\begin{abstract}
We introduce a geometric stencil selection algorithm for Laplacian in 3D that significantly improves octant-based 
selection considered earlier. The goal of the algorithm is to choose a small subset from a set of irregular points surrounding
a given point that admits an accurate numerical differentiation formula.
The subset serves as an influence set for the numerical approximation of the
Laplacian in meshless finite difference methods using either polynomial or kernel-based techniques.
Numerical experiments demonstrate a competitive performance of this method in comparison to the finite element method
and to other selection methods for solving the Dirichlet problems for 
the Poisson equation on several STL models. Discretization nodes for these domains are obtained either by 3D triangulations 
or from Cartesian grids or Halton quasi-random sequences.
 \begin{keywords}
 RBF-FD, meshless finite difference method, generalized finite differences
\end{keywords}
\end{abstract}

\section{Introduction}
\label{intro}

We consider the Dirichlet problem for the Poisson equation: Find a function $u: \overline{\Omega}  \to \mathbb{R}$ such that
\begin{equation}
	\begin{array}{l}
 	\Delta u = f\quad \text{in}\quad \Omega, \\
	\quad u = g\quad \text{on}\quad \partial \Omega. 
	\end{array}
\label{pr}
\end{equation}
where $\Delta$ is the Laplacian operator, $\Omega  \subset {\mathbb{R}^3}$  a domain, $f$  a real-valued function defined in $\Omega$, 
and $g$ a real-valued function defined on $\partial \Omega $. 

Meshless finite difference methods for \eqref{pr} are defined as follows.
Let $\Xi$ be  a
finite set of discretization nodes, where  $\Xi  \subset \overline \Omega$, and $\Xi_{\rm int}: = \Xi  \cap \Omega$,
$\partial \Xi : = \Xi  \cap \partial \Omega$. 
For each $\zeta \in \Xi_{\rm int}$, choose a \emph{set of influence}
$\Xi_{\zeta}:= \{\xi_1,\xi_2, \ldots,\xi_k\}\subset\Xi$ with $\xi_1 = \zeta$,
and approximate  $\Delta u\left( \zeta \right)$ by a numerical differentiation formula
\begin{equation} \label{delta_u}
\Delta u\left( \zeta  \right) \approx  \sum\limits_{\xi  \in {\Xi _\zeta }} {w_{\zeta ,\xi }} u\left( \xi  \right),
\quad \zeta  \in {\Xi_{\rm int}}.
\end{equation}
Thus, in contrast to the classical Finite Difference Method, \emph{stencils} defined by the sets of influence and
corresponding numerical differentiation weights $w_{\zeta ,\xi }$ are chosen individually for each node $\zeta$, rather than
obtained by a simple scaling of one and the same stencil. This removes the need for a grid structure in $\Xi$ and 
hence allows a meshless method.

The discretization of the problem \eqref{pr} is given by the square linear system
\begin{equation}  
\begin{array}{l}
\sum\limits_{\xi  \in {\Xi_\zeta }} {w_{\zeta ,\xi }} \hat{u}_\xi  = f\left( \zeta  \right),\quad \zeta  \in {\Xi _{\rm int}};\\
\quad \quad \quad \;\;\hat{u}_\xi  = g\left( \xi  \right),\quad \xi  \in \partial\Xi .
\end{array}
\label{prd}
\end{equation}
By solving the system \eqref{prd} we obtain a vector $[\hat{u}_\xi]_{\xi \in \Xi}$ that serves as 
an approximation of the solution $u$ of  \eqref{pr} at all nodes $\xi \in \Xi$.

A major motivation for meshless methods is that
no connection information between nodes in
$\Xi$ is needed, which simplifies the task of \emph{node generation} in comparison to
\emph{mesh generation}, where the nodes are connected into  grid structures or  elements that have to meet strict
quality criteria. Obtaining good meshes needed in mesh-based methods is 
 often notoriously difficult in applications, 
especially when $\Omega$ is defined by a complicated CAD model, or when the mesh has to be quickly adapted to moving boundaries or
singularities of the solution evolving in time.

Any particular meshless finite difference (mFD) method  is defined by specific  \emph{stencil selection algorithm}  
 for choosing the  sets of influence $\Xi_\zeta \subset \Xi$, $\zeta \in \Xi_{\rm int}$, as well as an algorithm for computing
 numerical differentiation weights $w_{\zeta,\xi}$, 
$\zeta\in\Xi_{\rm int}$, $\xi\in\Xi_{\zeta}$. The weights can be found by the traditional technique
of requiring exactness of the formula \eqref{delta_u} for polynomials of certain degree, as e.g.\ in \cite{Jensen72}, %
or by a radial basis function
method (RBF-FD), see e.g.\ \cite{FFprimer15} %
and references therein. The accuracy of numerical differentiation formulas on irregular points obtained by various methods has 
been investigated in \cite{DavySchaback16,DavySchabackMND,DavySchabackOPT,D21}.

The density of the linear system \eqref{prd} is determined by the size $k$ of the sets of influence $\Xi_{\zeta}$ 
since the corresponding row of the system matrix has at most $k$ nonzero entries. 
Since the density of the matrix is directly related to the cost of solving the resulting linear system, it is important 
to keep $k$ as small as possible, even if more computation time is spent on selection. Indeed, the selection is done 
independently for each node $\zeta$, and therefore can be done in parallel without any complications, whereas the linear
system may be very large for a large scale problem, so that higher sparsity is an important advantage.

In \cite{DavyOanh11, OanhDavyPhu17} we developed adaptive meshless RBF-FD %
algorithms for solving elliptic equations with point singularities in 2D, with sets of influence consisting of just $k=7$
 points, which delivers system matrices of about the same density as in the finite element method with linear shape
functions. The accuracy of the solution is also similar to the finite element counterpart, while our method is purely
meshless.

A general algorithm for the selection of sparse sets of influence based on the pivoted QR decomposition of the 
polynomial collocation matrices has been introduced in \cite{D19arxiv} and tested on elliptic interface problems
in \cite{DavySaf21}, showing convincing  performance especially in the context of a high order polynomial-based 
mFD method.

In the recent paper  \cite{DOT20} we investigated a simple method for selecting sets of influence in 3D 
that generalizes the well-known two-dimensional four-quadrant criterium of \cite{LiszOrk80}. By selecting 
two points in each Cartesian octant centered at $\zeta$, we obtain $\Xi_\zeta$ with $k=17$, which leads to
only slightly denser system matrices than those arising in the 3D finite element method with linear shape functions.

In this paper we suggest an improvement of the selection algorithm of \cite{DOT20} and compare 
its numerical performance with several alternative approaches for the selection of sets of influence in the meshless
finite difference method, as well as with the finite element method. Comparison is performed on Poisson equations with
Dirichlet boundary conditions for the unit ball, with known exact solution, and for three non-convex STL models with 
homogeneous Dirichlet boundary conditions, where we rely on reference solutions obtained with the finite element method
on fine triangulations. We consider several methods of node generation, including nodes obtained by mesh generation
(MATLAB PDE Toolbox triangulations; non-optimized Delaunay triangulations by Gmsh) as well as Cartesian grids and Halton
quasi-random points. The results indicate that the new method works well even on less regular nodes and competes 
successfully with the finite element method and alternative selection algorithms for mFD.

The paper is organized as follows. Section~\ref{improve} presents the new selection algorithm,
whereas Sections~\ref{weights} and \ref{sten} briefly describe weight computation and alternative selection
methods. Section~\ref{numexp} is devoted to numerical experiments, and Section~\ref{concl} provides a short conclusion.

\section{Improvement of octant-based selection}\label{improve}

The goal of a selection algorithm is to obtain suitable influence sets 
$\Xi_\zeta= \{\xi_1, \ldots,\xi_k\}\subset\Xi$,  $\xi_1 = \zeta$, with as small as possible $k$, 
without expectation that the node set $\Xi$ is highly regular.
The simplest selection algorithm builds $\Xi_\zeta$ by using $\zeta$ and its $k-1$ nearest nodes in $\Xi$.  
According to numerical experiments in \cite{BFFB17}, even on very nicely distributed nodes obtained 
by rather expensive node generation methods, the number of points in $\Xi_\zeta$ has to be is at least twice the dimension
of polynomials of degree $p$ in order to achieve convergence of the method \eqref{prd} with order $\mathcal{O}(h^p)$, 
where $h$ measures the spacing of the nodes. In 3D this approach requires $k\ge20$ if we pursue a method comparable 
in accuracy to the finite element method with linear shape functions whose convergence order is $\mathcal{O}(h^2)$. 

We consider a selection algorithm as successful if (a) it produces $\Xi_\zeta$ with $k$ less or equal 20 on cheaply generated 
nodes obtained without node improvement by repulsion and similar techniques that lead to highly regular nodes as in 
\cite{BFFB17}, and (b) it shows similar accuracy to the finite element method with the same number of degrees of freedom.
The octant-based selection algorithms of \cite{DOT20} have been successful in this sense on nodes produced as vertices of
a triangulation generated by MATLAB PDE Toolbox. However, these nodes are rather nicely distributed and therefore do not fully
demonstrate the advantages of a meshless method. As we will see in Section~\ref{numexp}, simple octant-based selection
of \cite{DOT20} often fails on less regular nodes. We suggest here an improvement of this method that works well on irregular
nodes obtained from non-optimized Delaunay triangulations or by combining Cartesian or Halton nodes inside a 3D domain with 
some discretization of its boundary.

A ``geometric'' selection method, such as our 2D algorithms in \cite{DavyOanh11, OanhDavyPhu17} seeks  
a compromise between two goals that the nodes $\xi_i$ are distributed as evenly as possible around $\zeta$ and 
in the same time as near to $\zeta$ as possible. 
In the earlier versions of our octant-based algorithms \cite{DOT20} either $n=2$ closest neighboring
nodes are chosen in each coordinate octant around the given node $\zeta\in\Xi$ \cite[Algorithm 1]{DOT20}, 
or only one closest node is taken in each of 16 half-octants \cite[Algorithm 2]{DOT20}, leading again to 2 nodes in each octant,
and $k\le 17$.
These algorithms are fast and their results are quite good for relatively nice discretizations $\Xi$.
However, the nodes selected by this method may still be poorly distributed around $\zeta$, or even contain clusters, while 
wide ranges of directions from $\zeta$ may receive no nodes in $\Xi_\zeta$ if $\Xi$ is rather irregular.

In the new algorithm, in order to improve the angular distribution of $\xi_i$ around $\zeta$, 
we first select a larger number $n>2$ of
candidate nodes in each of the eight octants. This produces a cloud of neighbors that better cover all
directions from $\zeta$, but are too many, for example up to 24 in the case $n=3$. Then, in the second stage of the algorithm, 
we go over these candidate nodes and choose at most $k-1$ of them to include in the final set of  influence $\Xi_\zeta$,
so that nodes that are better separated from already selected nodes, but also those  that are closer to $\zeta$ are preferred.

To control separation, we choose a \emph{standard distance $\rho _\zeta$} defined as 
$$
\rho _\zeta :=  \frac{\delta}{6}\,\sum_{i=1}^6\|\zeta-\xi_i\|,$$
where $\xi_1,\ldots,\xi_6$ are the six closest to $\zeta$ nodes in $\Xi\setminus\{\zeta\}$, 
$\left\|  \cdot  \right\|$ denotes the Euclidean norm in
${\mathbb{R}^3}$, and $\delta$ is a user-defined parameter satisfying $0<\delta<1$.

 Recall that an \emph{octant} centered at the origin in 3D is a cone that consists of points $(x_1,x_2,x_3)\in\RR^3$ such that
$\epsilon_i x_i\ge 0$, $i=1,2,3$, where $\epsilon_i\in\{-1,1\}$. The octants are not disjoint, therefore we reduce 
the boundary of some
octants to make sure that for each point with e.g.\ $x_1=0$ there is only one octant that contains it, see \cite{DOT20}.
The points in a neighborhood of $\zeta$ are classified according to 8 octants centered at $\zeta$.
In \cite[Algorithm 2]{DOT20} we split each octant in two \emph{half-octants}, leading to the partition of the space around
$\zeta$ in 16 cones, and in this paper we also consider 
\emph{one-third-octants} defined by the condition $\max\{|x_1|, |x_2|, |x_3|\}=|x_i|$, $i=1,2,3$, with appropriate reduction
of the boundaries of some of them, which partitions the space into 24 disjoint cones.

We now describe the new selection algorithm that will be referred to as \texttt{oct-dist} in the rest of the paper.

\begin{algorithm}[\texttt{oct-dist}]\label{oct-dist}\rm
\textit{Input}: $\Xi$ and $\zeta\in \Xi_{\rm int}$.
\textit{Output}: $\Xi_{\zeta}$.\\
\textit{Parameters}: 
\begin{itemize}
\item %
$m$: the number of nodes in the initial local cloud including $\zeta$; 
\item   %
$k$: the target number of nodes in $\Xi_\zeta$;
\item $s \in \{1, 2, 3\}$: the number of subdivisions of octants;
\item   %
$n$: the number of candidate nodes in each octant, $n=s\nu$ for some $\nu\in\NN$;
\item  $0<\delta<1$: the standard distance tolerance.  
\end{itemize}

\begin{enumerate}[\hspace*{-10pt}I.]
\item
\label{preat} {Choose a set $\Xi_{\rm cand}\subset \Xi \setminus \{\zeta\}$ of candidate neighbors.}

\begin{enumerate}[1.]
\makeatletter
\renewcommand*\p@enumii{\theenumi.}
\makeatother

\item\label{preat-I-1}
Choose the initial cloud $\Xi_{\rm Init}:= \{\xi_1, \xi_2, \ldots , \xi_{m-1}\} \subset \Xi \setminus \{\zeta\} $ 
consisting of $m-1$ nodes closest to $\zeta$, sorted by increasing distance to $\zeta$. 
\item 
Compute $\rho _\zeta :=  \frac{\delta}{6}\,\sum_{i=1}^6\|\zeta-\xi_i\|$.

 \item \label{preat-I-4}
Determine  $8s$ sets $\tilde {\rm O}_j$, $j=1, 2,\ldots, 8s$, corresponding to the octants,
for $s=1$, half-octants  for $s=2$, or one-third-octants  for $s=3$, by collecting in each $\tilde {\rm O}_j$ at most $\nu$ 
nodes in $\Xi_{\rm Init}$ closest to $\zeta$ and lying in the corresponding octant,  half-octant or one-third-octant.
We set $\Xi_{\rm cand} :=\bigcup_{j=1}^{8s}\tilde {\rm O}_j$.

 \item \label{preat-I-5}
Determine subsets ${\rm O}_j = \{ \xi^j_1, \xi^j_2, \dots\}$, $j=1, 2,\ldots, 8$, of $\Xi_{\rm cand}$ consisting
of the nodes lying in the eight octants as follows. If $s=1$, then ${\rm O}_j =\tilde {\rm O}_j$; otherwise ${\rm O}_j$ is the union
of $s$ sets of the type $\tilde {\rm O}_i$. Assume that the numbering of the nodes is such that
$\| \zeta -\xi^j_1  \| \leq\| \zeta -\xi^j_2\| \leq \cdots\leq \| \zeta - \xi^j_{n_j}\|$, with $0\le n_j\le n$. 
Note that some ${\rm O}_j$ may be empty.

\item \label{preat-return}
\emph{If} $\#\Xi_{\rm cand} \le k-1$,  \emph{then} STOP and return $\Xi_\zeta : = \Xi_{\rm cand}  \cup\{ \zeta\}$.

\end{enumerate}

\item\label{deter} {Choose the influence set $\Xi_\zeta\subset \Xi_{\rm cand}\cup\{\zeta\}$.}

\emph{Initialization}: $\Xi_\zeta := \emptyset$.%
	
\begin{enumerate}[1.]
\makeatletter
\renewcommand*\p@enumii{\theenumi.}
\makeatother

\item \label{deter-1}
\emph{For} $j=1$ to $8$: \\
\emph{If} ${\rm O}_j \ne\emptyset$, \emph{then} $\Xi_\zeta : = \Xi_\zeta \cup \{\xi^j_1\}$.

 \item \label{deter-2}
Set
$\Xi_{\rm cand} :=\Xi_{\rm cand}  \setminus \Xi_\zeta $, such that $\Xi_{\rm cand}=\{ \bar\xi_1, \bar\xi_2, \dots\}$, where 
$\bar \xi_1, \bar \xi_2,\dots$ are sorted by increasing distance to $\zeta$. 

\item\label{deter-3}
\emph{For} $i=1$ to $\# \Xi_{\rm cand}$:\\
\emph{If} $ \dist( \bar\xi_i,\Xi_\zeta) \geq \rho _\zeta $, \emph{then:}
\begin{enumerate}[a.]
\item 
 $\Xi_\zeta : = \Xi_\zeta \cup \{ \bar\xi_i\}$.
\item \label{deter-return}
 \emph{If} $\#\Xi_\zeta = k-1$,  \emph{then} STOP and return $\Xi_\zeta : = \Xi_\zeta  \cup\{ \zeta\}$.
\end{enumerate}

\item  \label{deter-3i}
Set 
$\rm \rho _\zeta:=  \delta \rho _\zeta$ and GOTO \ref{deter-2}. 

\end{enumerate}
 
\end{enumerate}
\end{algorithm}
 
\subsubsection*{Remarks}

\begin{enumerate}[1.]
\item  Based on our numerical experiments, $m = 100$ seems a good choice for this parameter;
$s=1$ with $n=3$ works well in most cases, but for less regularly distributed nodes larger
values for $s$ and $n$ are sometimes advantageous. We suggest $k=17$ as default value, with
smaller $k$ sometimes giving better results on more regular nodes. We choose
$\delta$ between 0.7 and 0.9, with a higher value for more complicated 3D shapes, 
see more details in Section~\ref{numexp}.

\item 
In the first stage of the algorithm the initial cloud of $m-1$ nearest neighbors of $\zeta$ is reduced to  the set 
$\Xi_{\rm cand}$ of at most $8n$ nearby nodes 
that represent the directions of all octants from $\zeta$.  Using a larger parameter $s$ makes these candidate neighbors more evenly distributed at the
expense of potentially larger distances to $\zeta$. If $\Xi_{\rm cand}$ contains less than $k$ nodes, then we skip the second
stage and return  $\Xi_\zeta : = \Xi_{\rm cand}  \cup\{ \zeta\}$. Note that we usually choose $n\ge k/8$, so that this
situation may only occur if some octants contain too few nodes of $\Xi$ in the vicinity of $\zeta$.

\item
In the second stage at most $k-1$ candidate nodes of $\Xi_{\rm cand}$ are selected into the set of influence $\Xi_\zeta$. 
First, the closest node in each octant is selected and removed from  $\Xi_{\rm cand}$. 
In Step~\ref{deter-3} we go over the remaining candidate nodes according to their distance to $\zeta$, and add them 
to $\Xi_\zeta$ if they are separated from already selected nodes by at least the standard distance $\rho_\zeta$.
This loop is terminated if $k-1$ candidate notes are selected, in which case we add $\zeta$ itself to $\Xi_\zeta$ and return
this set. Otherwise, we relax the separation restriction by setting $\rm \rho _\zeta:=  \delta \rho _\zeta$, and repeat the
selection loop of Step~\ref{deter-3} with the remaining candidate nodes. This process will eventually
terminate at Step~\ref{deter-3}.b because the number of
candidate nodes at the second stage is at least $k$, and the distance condition gradually weakens.

\item In the earlier versions of octant based selection \cite[Algorithms 1 and 2]{DOT20}, we modified 
the sets $\Xi_\zeta$ in the following way when  $\Omega$ is not convex: For
each `invisible' $\xi\in\Xi_\zeta\setminus\{\zeta\}$  such that the segment $(\zeta,\xi)$ intersects the boundary 
$\partial\Omega$ of $\Omega$, we replaced $\xi$  by the point $\xi'\in (\zeta,\xi)\cap\partial\Omega$ closest to
$\zeta$. We do not apply this step here because the new points $\xi'$ tend to significantly reduce the separation
distance of the nodes in $\Xi_\zeta$ and often lead to stencils of poor quality. Clearly, using nodes $\xi$ with 
$(\zeta,\xi)\cap\partial\Omega\ne\emptyset$ may be harmful on domains that possess thin holes or cracks, 
with big jumps in the solution across them. A special treatment may be needed in this case, for example the initial 
cloud $\Xi_{\rm Init}$ at Step~\ref{preat-I-1} may be chosen in such a way that it does contain those neighbors 
$\xi$ for which $\zeta$ and $\xi$ are not connected by a short path within $\Omega$. 
However, `invisible' nodes do not normally seem to be 
a problem if the discretization set $\Xi$ is sufficiently dense in $\Omega$.

\end{enumerate}

\section{RBF-FD weights}\label{weights}
Once a set of influence $\Xi_\zeta$ is selected, the weights $w_{\zeta,\xi}$ of the numerical differentiation formula
\eqref{delta_u} can be computed by the radial basis function method (RBF-FD). 

Let $\Phi :{\mathbb{R}^3} \to \mathbb{R}$ be a positive definite or conditionally positive definite radial basis function,
$\Phi \left( x \right) =\varphi \left( {\left\| x \right\|} \right)$, $x \in {\mathbb{R}^3}$, with continuous
$\varphi :\mathbb{R}_{+}\rightarrow \mathbb{R}$,  
see \cite{Buhmann03,Fasshauer07,Wendland} for further details on these functions.

For each $\zeta  \in {\Xi_{\rm int}}$, the weights $w_{\zeta,\xi}$ are determined by the exactness condition
\begin{equation}\label{exa}
\Delta s\left( \zeta  \right) =  \sum\limits_{\xi  \in {\Xi_\zeta }} {w_{\zeta ,\xi }} s\left( \xi  \right),
\end{equation}
required for all functions of the form
$$
s(x) =\sum_{\xi  \in {\Xi_\zeta }}a_{\xi}\Phi (x-\xi)+\sum_{j=1}^L c_jp_j(x), \quad a_\xi,c_j\in\RR,$$
satisfying the side condition
$$
\sum_{\xi  \in {\Xi_\zeta }} a_\xi p_i(\xi)=0,\quad i=1,\ldots,L,$$
where $\{p_1,\ldots,p_L\}$, $L={\ell+2\choose 3}$, is a basis for the linear space of trivariate polynomials of 
order  at most $\ell$ (that is, total degree at most $\ell-1$). As shown in \cite{D21}, the RBF-FD weights  $w_{\zeta,\xi}$ are uniquely determined 
by these conditions as soon as $\Phi$ is positive definite or conditionally positive definite of order at most $\ell$,
and there is any solution $\{v_{\zeta ,\xi }:\xi\in \Xi_\zeta\}$ for the polynomial exactness condition
 \begin{equation}\label{exap}
\Delta p_i\left( \zeta  \right) =  \sum\limits_{\xi  \in {\Xi_\zeta }} {v_{\zeta ,\xi }} p_i\left( \xi  \right),
\quad i=1,\ldots,L.
\end{equation}
We refer to \cite{D21} for a discussion of computational methods for  $w_{\zeta,\xi}$. 

As demonstrated in \cite{BFFB17}, particularly good results are obtained with polyharmonic RBF defined
by $\varphi(r)=r^\alpha$ with a polynomial term. In this paper we are interested in methods with expected convergence order
$\mathcal{O}(h^2)$, and therefore use  $\varphi(r)=r^5$ with quadratic polynomial term ($\ell=3$), which seem 
the most appropriate setting for this goal.

\section{Other stencil selection methods}\label{sten}

In addition to Algorithm~\ref{oct-dist} (\texttt{oct-dist}),
we consider for comparison several other stencil selection methods that produce influence sets of size not exceeding 20.
This number is motivated by the recommendation of \cite{BFFB17} that the  number of nodes in $\Xi_\zeta$ should be at least
the double of the polynomial dimension $L$, which is 10 in our experiments with the polyhamonic RBF $\varphi(r)=r^5$
and $\ell=3$.

\begin{description}

\item[\texttt{tet}] In this case a conforming triangulation $\triangle$ of $\Omega$ into tetrahedra with vertices at all nodes of
$\Xi$ is required, and $\Xi_\zeta$ consists of $\zeta$ and all nodes connected to it by an edge of  $\triangle$. 
This leads to a system matrix with sparsity pattern identical with that of the finite element method with 
linear shape functions on the same triangulation.  Note that although this approach does not seem particularly useful 
in the meshless context as it relies on a tetrahedral mesh, it is possible to obtain similar $\Xi_\zeta$ by a
local Delaunay triangulation in the neighborhood of $\zeta$, which may circumvent the main hurdles of mesh generation. 
This method produces sets of influence containing on average about 16 nodes in our experiments.

\item[\texttt{oct}] Stencil selection according to \cite[Algorithm 2]{DOT20}, where 100 nearest neighbors of $\zeta$ 
are sorted into 16 cones obtained by bisecting eight coordinate octants, and, in addition to $\zeta$, the closest node 
in each half-octant is taken into $\Xi_\zeta$. Clearly, the number of nodes in $\Xi_\zeta$ is at most 17.

\item[\texttt{20near}] The set of influence $\Xi_\zeta$ consists of 20 nearest neighbors of $\zeta$, including  
$\zeta$ itself. This is the cheapest and most commonly used method. 
As shown in \cite{DOT20} smaller number of neighbors leads to significantly larger errors. Occasionally we check the
performance of the methods \texttt{30near} and \texttt{40near} that select 30, respectively 40, nearest neighbors of $\zeta$.

\item[\texttt{pQR4sel}] Starting with  100 nearest neighbors of $\zeta$, we select $\Xi_\zeta$ by the pQR method of order 4,
as described in \cite{D19arxiv,DavySaf21}. This results in at most 20 nodes unless the polynomial exactness
condition \eqref{exap} of order 4 cannot be satisfied on 100 neighbors, in which case the method fails.

\end{description}

The above methods are applied only for selecting the sets of influence, while the weights $w_{\zeta,\xi}$ 
of \eqref{prd} are obtained by  the RBF-FD method described in Section~\ref{weights}, with $\varphi(r)=r^5$ and $\ell=3$.
On the other hand, pQR selection method introduced in \cite{D19arxiv} already generates weights satisfying the polynomial exactness condition \eqref{exap}, and 
corresponding meshless finite difference method performs well in the numerical experiments presented in 
\cite{DavySaf21}. Moreover, since the number of selected nodes in this case does not exceed the polynomial dimension $L$,
RBF-FD weights computed with the polynomial term of order $\ell$ are only rarely different from the pQR weights of the
same order. Hence, unless we choose a smaller order for RBF-FD, as in \texttt{pQR4sel}, there is no point in 
replacing pQR weights by almost always identical RBF-FD weights. 
Because of this we also consider the following two selection methods with pQR weights whose computation does not rely on
radial basis functions.

\begin{description}
\item[\texttt{pQR3}, \texttt{pQR4}] Starting with  100 nearest neighbors of $\zeta$, we select $\Xi_\zeta$ by the 
pQR method \cite{D19arxiv,DavySaf21} of order 3 or 4. Note that even if \texttt{pQR4} enforces exactness of the 
numerical differentiation weights
for polynomials of order $\ell=4$ rather than 3, this does not seem to lead to a higher order method, as seen in
numerical experiments below and in \cite{DavySaf21}. 
The number of nodes in $\Xi_\zeta$ is at most 10 for \texttt{pQR3} and at most
20 for \texttt{pQR4}.
\end{description}

\section{Numerical Experiments}\label{numexp}

In this section we present a number of numerical tests to compare the performance of the improved octant-based selection
method \texttt{oct-dist} of Algorithm~\ref{oct-dist} with other selection approaches for the meshless finite difference  
method, and with the finite element method. We use MATLAB 2017b in all experiments. 
In particular, finite element solutions are computed by PDE Toolbox \cite{PDEtool}.

We  check the performance of mFD methods on several different types of discretization nodes $\Xi$, 
produced either as the vertices of 3D triangulations, or as Cartesian grids or Halton quasi-random points
inside $\Omega$ combined with some discretization of the boundary.

In particular, for all test problems we run all methods
on sets $\Xi$ obtained as nodes of \emph{optimized triangulations} produced  by the command 
\texttt{generateMesh} of MATLAB PDE Toolbox,
with target maximum edge size $\texttt{Hmax}=H_02^{-i/3}$, $i=0,1,\ldots$,
where $H_0$ is chosen individually for each test problems because of the differences in the size and complexity
of the domains. The factor $2^{-1/3}$ is chosen such that the number of interior vertices approximately
doubles for the next value of \texttt{Hmax}, possibly after some warm-up on coarse triangulations. 
The other parameters of \texttt{generateMesh} are always set to $\texttt{Hmin}=\texttt{Hmax}/3$ and
\texttt{Hgrad}=1.5. We set \texttt{GeometricOrder} to 1 for the computation of the first order finite element 
solution \texttt{fem1} based on linear shape functions and all versions of the meshless finite difference method. 
It is set to 2 when we compute the second order finite element solution \texttt{fem2} that employs quadratic shape 
function.

In addition, for each test problem we also generate $\Xi_{\rm int}$ as uniform Cartesian \emph{grids} in the interior of 
$\Omega$ with target spacing $h=0.9\,\texttt{Hmax}$ and distance from the boundary at least $0.25h$. Similarly, as another
option, the nodes $\Xi_{\rm int}$ with the same target average spacing $h$ are generated from the 3D \emph{Halton}  quasi-random stream with the help of MATLAB command
\texttt{qrandstream}, scaled to a cube containing $\Omega$, and cleared of  points outside of $\Omega$ and those 
at a distance less than $0.25h$ from the boundary. In this cases we produce boundary nodes either by 
\emph{orthogonal projection} to the boundary of interior nodes at distance less than $h$, or by using boundary nodes of the
above optimized triangulation. Note that orthogonal projection  was successfully applied for boundary and interface 
discretization  in numerical experiments for elliptic interface
problems in \cite{DavySaf21}, even for higher order mFD methods with pQR selection.
When $\Omega$ is convex (Test Problem~\ref{pro11}) we use projected boundary nodes and also test the first order 
finite element method \texttt{fem1} on the same set $\Xi$. A mesh required for this method is a 3D Delaunay 
triangulation of the nodes in $\Xi$ produced  by MATLAB's command 
\texttt{delaunayTriangulation}. This provides us with an opportunity to see how various methods perform on vertices of an
\emph{unoptimized triangulation}. 

For non-convex domains \texttt{delaunayTriangulation} does not generate a 
triangulation of $\Omega$. In order to obtain unoptimized triangulations for such domains 
(Test Problems~\ref{pro19}--\ref{pro13}), we use Gmsh \cite{Gmsh}, where we switch off mesh smoothing 
and optimization by setting the parameters \texttt{-smooth} and \texttt{Mesh.Optimize} to zero. 
Similar to MATLAB mesh generation, we
produce a sequence of meshes by setting the maximum characteristic length parameter  to
$\texttt{clmax}=H_02^{-i/3}$, $i=0,1,\ldots$, with individual $H_0$ reported below for each test problem.

Once a numerical solution $\hat u=[\hat u_\xi]_{\xi\in\Xi}$ is computed by either meshless finite difference method 
\eqref{prd}, or by the finite element method, we measure the  accuracy of this solution by the 
relative root mean square (RRMS) error on the interior nodes, given by
\begin{equation}\label{rrms}
E_{\rm ref}= {\rm RRMS}(\hat u^{\rm ref},\hat u,\Xi_{\rm int}) 
:=\frac{\Big(\sum\limits_{\zeta  \in \Xi_{\rm int}}( \hat u^{\rm ref}_\zeta
- \hat u_\zeta)^2\Big)^{1/2}}{\Big(\sum\limits_{\zeta\in \Xi _{\rm int}}( \hat u^{\rm ref}_\zeta)^2\Big)^{1/2}}.
\end{equation}
where  the vector $\hat u^{\rm ref}=[\hat u^{\rm ref}_\xi]_{\xi\in\Xi}$ represents a reference solution, 
obtained by  either evaluating on $\Xi _{\rm int}$ the  exact solution $u$ of
\eqref{pr},  $\hat u^{\rm ref}_\xi:=u(\xi)$, $\xi\in\Xi$, if it is known analytically, 
or by interpolating to $\Xi _{\rm int}$ a numerical solution from a much finer set of nodes. For this we employ 
MATLAB's command \texttt{scatteredInterpolant} with default settings, which makes use of piecewise linear interpolation 
over a triangulation of the nodes.
Note that we write $E_{\rm ref}=\texttt{NaN}$ in the tables below when the weights of Section~\ref{weights}
cannot be found for some $\zeta\in\Xi_{\rm int}$, which only happens when there are no weights satisfying the 
polynomial exactness condition \eqref{exap}. On rare occasions we write $E_{\rm ref}=\texttt{Inf}$, which means that 
MATLAB failed to solve the system \eqref{prd} because its matrix was singular to working precision.

In addition to the errors, we provide information about  the density of the system 
matrices since various selection methods lead to sets of influence of different sizes.  The density of a matrix $A \in \mathbb{R}^{n \times n}$ is the average number of nonzero entries
per row, computed as 
\begin{equation}\label{density}
\texttt{density}=\texttt{nnz}(A)/n,
\end{equation}
 where {\tt nnz}$(A)$ is the total number of nonzero entries in $A$.
Note that we measure the density of the matrix of \eqref{prd} and that of the finite element stiffness matrix after eliminating the Dirichlet
boundary conditions. Therefore the density is lower for coarse discretization where a more significant percentage
of the nodes is located on the boundary.

In most cases the linear system \eqref{prd} is solved by a direct method built into MATLAB's backslash command. 
However, for Test Problem~\ref{pro19} we also run an iterative solver in order to see how the number of iterations 
depends on the choice of the selection method.

\begin{testproblem}\label{pro11}%
Poisson equation $\Delta u = 3 e^{x+y+z}$ in the unit ball  
 $\Omega = \{(x,y,z)\in {\mathbb{R}^3}: x^2 + y^2 + z^2  < 1 \}$  
 with Dirichlet boundary conditions chosen such that the exact solution is $u(x, y,z) = e^{x+y+z}$.
\end{testproblem}

In the first experiment we obtain $\Xi$ as vertices of the optimized triangulations produced  by 
\texttt{generateMesh} as described above, with $H_0=0.25$. The results are presented in Table~\ref{pr11fem}.
The column marked \texttt{fem1} shows the error $E_{\rm ref}$ for the first order finite element method, whereas the other
columns correspond to the meshless finite difference method with various versions of selection algorithms, as described in
Sections~\ref{improve} and \ref{sten}. We reserved two columns for Algorithm~\ref{oct-dist},  \texttt{oct-dist13} and 
\texttt{oct-dist17}, depending on the value of the parameter $k=13$ or 17. Other parameters of Algorithm~\ref{oct-dist}
are the same in both cases: $s=1$, $n=3$, $\delta=0.9$. The last row of the table gives the density of the 
system matrix on the finest set of nodes.

In addition, Table~\ref{pr11fem2} presents the errors and density for 
the second order finite element method on the nodes of the optimized triangulations obtained with 
$H_0=0.5$ and $\texttt{GeometricOrder}=2$.

\begin{table}[htbp!]\scriptsize %
\centering
\renewcommand{\arraystretch}{1.25}
\begin{tabular}{|c||c||c|c|c|c|c|c||c|c|}
\hline
\multirow{2}{.7cm}{$\# {\Xi _{\rm int}}$} & FEM   
 &	\multicolumn{6}{|c||}{Polyharmonic $r^5$ with quadratic polynomial} & \multicolumn{2}{c|}{Pivoted QR}\\
\cline{2-10}
		& {\scriptsize fem1}	& {\scriptsize tet} &	{\scriptsize oct} & {\scriptsize oct-dist13} & {\scriptsize oct-dist17}
		& {\scriptsize 20near} & {\scriptsize pQR4sel} & {\scriptsize pQR3} & {\scriptsize pQR4}\\
\hline
245	&	5.4e{\tt-}03	&	2.6e{\tt-}03	&	2.9e{\tt-}03	&	2.6e{\tt-}03	&	2.7e{\tt-}03	&	2.8e{\tt-}03	&	3.1e{\tt-}03	&	3.2e{\tt-}03	&	2.9e{\tt-}03\\
567	&	3.3e{\tt-}03	&	NaN	&	2.1e{\tt-}03	&	1.6e{\tt-}03	&	1.8e{\tt-}03	&	1.7e{\tt-}03	&	1.3e{\tt-}03	&	1.9e{\tt-}03	&	1.7e{\tt-}03\\
1142	&	1.8e{\tt-}03	&	8.7e{\tt-}04	&	5.8e{\tt-}04	&	9.5e{\tt-}04	&	8.0e{\tt-}04	&	1.0e{\tt-}03	&	6.3e{\tt-}04	&	9.4e{\tt-}04	&	1.0e{\tt-}03\\
2523	&	1.0e{\tt-}03	&	5.3e{\tt-}04	&	4.2e{\tt-}04	&	6.6e{\tt-}04	&	5.5e{\tt-}04	&	5.7e{\tt-}04	&	3.5e{\tt-}04	&	7.8e{\tt-}04	&	5.9e{\tt-}04\\
5207	&	5.9e{\tt-}04	&	3.1e{\tt-}04	&	2.6e{\tt-}04	&	4.0e{\tt-}04	&	2.9e{\tt-}04	&	3.1e{\tt-}04	&	2.1e{\tt-}04	&	3.2e{\tt-}04	&	3.5e{\tt-}04\\
10780	&	3.8e{\tt-}04	&	1.8e{\tt-}04	&	1.5e{\tt-}04	&	1.8e{\tt-}04	&	9.7e{\tt-}05	&	1.4e{\tt-}04	&	1.0e{\tt-}04	&	2.1e{\tt-}04	&	2.2e{\tt-}04\\
21730	&	2.3e{\tt-}04	&	NaN	&	8.8e{\tt-}05	&	1.2e{\tt-}04	&	6.8e{\tt-}05	&	9.0e{\tt-}05	&	5.9e{\tt-}05	&	1.1e{\tt-}04	&	1.4e{\tt-}04\\
43956	&	1.4e{\tt-}04	&	NaN	&	4.5e{\tt-}05	&	7.2e{\tt-}05	&	3.6e{\tt-}05	&	4.1e{\tt-}05	&	2.7e{\tt-}05	&	8.1e{\tt-}05	&	8.5e{\tt-}05\\
88936	&	8.9e{\tt-}05	&	NaN	&	3.0e{\tt-}05	&	4.4e{\tt-}05	&	2.2e{\tt-}05	&	2.7e{\tt-}05	&	1.3e{\tt-}05	&	5.4e{\tt-}05	&	5.3e{\tt-}05\\
\hline
density:	&	14.6	&	14.6	&	16.5	&	12.7	&	16.5	&	19.4	&	19.4	&	9.8	&	19.4\\
\hline
\end{tabular}
\caption{Test Problem~\ref{pro11}: RRMS errors $E_{\rm ref}$ on interior nodes of the optimized
triangulations with $H_0=0.25$.
The last row shows the density of the system matrix on the finest set of nodes. 
}
\label{pr11fem}
\end{table}

\begin{table}[htbp!]\scriptsize %
\centering
\renewcommand{\arraystretch}{1.25}
\begin{tabular}{|c|c|c|c|c|c|c|c|c|c||c|}
\hline
$\# {\Xi _{\rm int}}$ & 272	&	824	&	1113	&	2357	&	5127	&	10138	&	21821	&	44238	&	90396	&	density\\
\hline
fem2	&	9.4e{\tt-}04	&	3.8e{\tt-}04	&	2.8e{\tt-}04	&	1.2e{\tt-}04	&	6.1e{\tt-}05	&	2.8e{\tt-}05	&	1.3e{\tt-}05	&	5.8e{\tt-}06	&	3.0e{\tt-}06	&	27.3\\
\hline
\end{tabular}
\caption{Test Problem~\ref{pro11}: RRMS error $E_{\rm ref}$ of the quadratic finite element method 
for the optimized triangulations with $H_0=0.5$. 
The last column shows system matrix density on the finest set of nodes.}
\label{pr11fem2}
\end{table}

The errors for all versions of the meshless finite difference method, including \texttt{tet} selection whenever this method does not
fail, are smaller (sometimes more than 6 times smaller)
than those of the finite element method \texttt{fem1} on the same nodes.  

The errors of \texttt{fem2} in Table~\ref{pr11fem2} on comparable sizes of $\Xi _{\rm int}$
are significantly better than for any method of Table~\ref{pr11fem}. 
This is explained by the fact that the solution $u$ is a very smooth, infinitely
differentiable function, so that higher order methods such as \texttt{fem2} should be beneficial.  However, we also notice
that the density 27.3 of the system matrix of \texttt{fem2}  is significantly higher than the densities seen in 
Table~\ref{pr11fem}.

\begin{figure}[htbp!]
\begin{center}
\vspace*{-1cm}
\subfigure[Optimized triangulation]{\includegraphics[width=9.6cm,height=6.4cm]{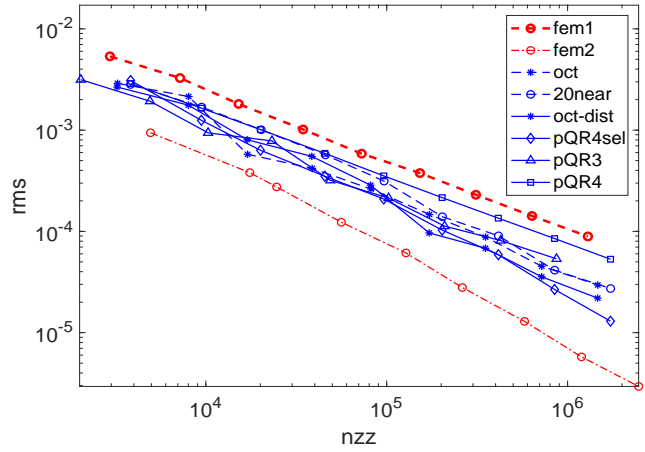}}
\subfigure[Interior grid and projected boundary nodes]{\includegraphics[width=9.6cm,height=6.4cm]{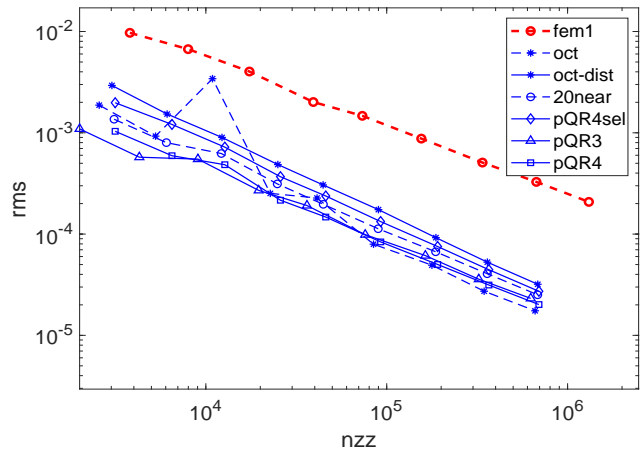}}
\subfigure[Halton interior and projected boundary nodes]{\includegraphics[width=9.6cm,height=6.4cm]{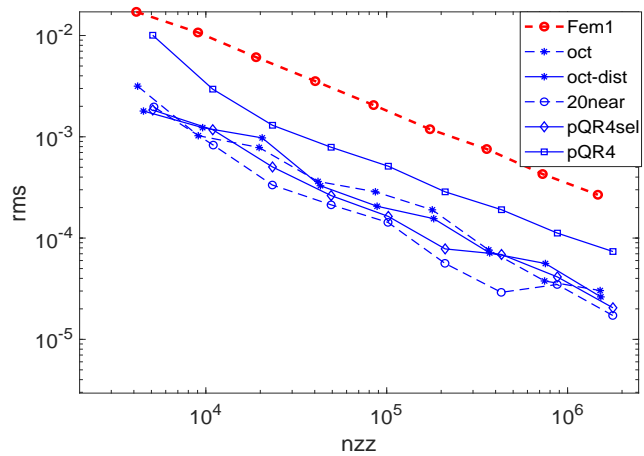}}
\caption{Test Problem~\ref{pro11}: 
RRMS errors as functions of nominal \texttt{nnz}, 
corresponding to (a)~Tables~\ref{pr11fem} and \ref{pr11fem2}, where
\texttt{oct-dist} stands for \texttt{oct-dist17}, (b) Table~\ref{pr11grid}, (c) Table~\ref{pr11halton}
}
\label{Pr1_tables}
\end{center}
\end{figure}

We notice that smaller errors in Table~\ref{pr11fem} clearly correlate with higher density. 
Comparing methods of the same density, we see that \texttt{oct-dist17} significantly outperforms \texttt{oct},
and \texttt{pQR4sel}  outperforms \texttt{20near} and  \texttt{pQR4}. In order to take into account the
differences in density,  we compare in Figure~\ref{Pr1_tables}(a) the errors of different methods as functions of the nominal 
number of nonzeros (nominal \texttt{nnz}) of the
system matrix, computed as the number of interior nodes $\# {\Xi _{\rm int}}$ multiplied by the density on the finest 
set of nodes (provided in the tables). Note that the number of nonzeros in the
system matrix determines the cost of matrix-vector multiplication, and therefore the cost per iteration when solving the
system by iterative methods. Figure~\ref{Pr1_tables}(a) does not include \texttt{tet} because it fails on several sets $\Xi$,
and for \texttt{oct-dist} we included the better performing version with $k=17$ and left out the one with $k=13$.

We see that the best performing selection methods for the nodes of the optimized triangulation are \texttt{oct-dist} with $k=17$ and \texttt{pQR4sel}.
in particular, both of them produce significantly more accurate solutions than \texttt{fem1}.

In the next experiment we generate $\Xi$ as a uniform grid in the interior of $\Omega$, with boundary nodes produced
by orthogonal projection as explained in the beginning of this section. 
The results are presented in Table~\ref{pr11grid} and Figure~\ref{Pr1_tables}(b). 
For the meshless finite difference method we take advantage of the grid structure of the nodes,
 and use for the Laplacian in \eqref{delta_u} the classical 7-node stencil,
\begin{align}\label{7star}
\begin{split}
&\Xi^{\rm 7star}_\zeta=\{\zeta,\zeta\pm (h,0,0),\zeta\pm (0,h,0),\zeta\pm (0,0,h)\},\\
&w_{\zeta,\zeta}= -6h^{-2},\quad w_{\zeta,\xi}=h^{-2},\;\xi\in \Xi_\zeta\setminus\{\zeta\},
\end{split}
\end{align}
whenever possible.
These weights are exact for polynomials of order 4. The same weights are also obtained by the radial basis exactness condition
\eqref{exa} on $\Xi^{\rm 7star}_\zeta$ when $\ell=3$ or 4.  Thus, specific stencil selection and weight computation
methods described in Sections~\ref{improve}--\ref{sten} are only applied for $\zeta\in \Xi_{\rm int}$ near the boundary such that
 $\Xi^{\rm 7star}_\zeta\not\subset\Xi_{\rm int}$. Therefore the density of the system matrix on finer discretizations 
 only slightly exceeds 7.
The results for the finite element method \texttt{fem1} are obtained on a 3D Delaunay triangulation of the nodes in $\Xi$
created by MATLAB command \texttt{delaunayTriangulation}, which produces a triangulation of $\Omega$ since 
this domain (the unit ball) is convex. 
We use $k=18$ for \texttt{oct-dist}, with values of other parameters of this method the same as above.
The errors demonstrated by all versions of the mFD method, except of  \texttt{tet} that again fails on
several sets of nodes, are very close to each other and significantly better than the errors of the finite element method,
especially if the density of the system matrices is taken into account. In Table~\ref{pr11grid} we also included selection of 30 nearest nodes
(\texttt{30near}), which does not increase the density of the system matrix essentially because of the 7-node stencils
still used in most places. The errors are slightly better than those for \texttt{20near}. Therefore larger sets of influence
may be recommended in this case. However, this improvement may be attributed to the high smoothness of $u$ as it
does not happen for more interesting problems considered below, as will be demonstrated for Test Problem~\ref{pro19}.
Comparison of Figures~\ref{Pr1_tables}(a) and \ref{Pr1_tables}(b) shows that the error to \texttt{nnz} ratio is better for mFD methods on
interior grid nodes than for those on the optimized triangulation, whereas \texttt{fem1} shows the opposite behavior. Note
that the error of the finite element method is significantly higher in comparison to the optimized triangulation, and the
density of the system matrix is not decreased on the grid nodes, whereas mFD methods deliver similar accuracy on both types
of nodes, but the density is much lower in the gridded case.

\begin{table}[htbp!]\scriptsize
\centering
\renewcommand{\arraystretch}{1.25}
\begin{tabular}{|c||c||c|c|c|c|c|c||c|c|}
\hline
\multirow{2}{.7cm}{$\# {\Xi _{\rm int}}$} & FEM   
 &	\multicolumn{6}{|c||}{Polyharmonic $r^5$ with quadratic polynomial} & \multicolumn{2}{c|}{Pivoted QR}\\
\cline{2-10}
		& {\scriptsize fem1} & {\scriptsize tet} 	&	{\scriptsize oct} & {\scriptsize oct-dist} 
		& {\scriptsize 20near} & {\scriptsize 30near} & {\scriptsize pQR4sel} & {\scriptsize pQR3} & {\scriptsize pQR4}\\
\hline
304	&	9.7e{\tt-}03	&	9.9e{\tt-}04	&	1.9e{\tt-}03	&	2.9e{\tt-}03	&	1.4e{\tt-}03   &	6.6e{\tt-}04	&	2.0e{\tt-}03	&	1.1e{\tt-}03	&	1.0e{\tt-}03\\
624	&	6.7e{\tt-}03	&	6.4e{\tt-}04	&	9.3e{\tt-}04	&	1.5e{\tt-}03	&	8.0e{\tt-}04   &	4.1e{\tt-}04	&	1.2e{\tt-}03	&	5.7e{\tt-}04	&	5.9e{\tt-}04\\
1308	&	4.0e{\tt-}03	&	3.5e{\tt-}04	&	3.4e{\tt-}03	&	9.0e{\tt-}04	&	6.2e{\tt-}04 &	3.9e{\tt-}04	&	7.2e{\tt-}04	&	5.5e{\tt-}04	&	4.8e{\tt-}04\\
2822	&	2.0e{\tt-}03	&	2.1e{\tt-}04	&	2.5e{\tt-}04	&	4.9e{\tt-}04	&	3.1e{\tt-}04 &	2.0e{\tt-}04	&	3.7e{\tt-}04	&	2.7e{\tt-}04	&	2.2e{\tt-}04\\
5196	&	1.5e{\tt-}03	&	1.4e{\tt-}04	&	2.3e{\tt-}04	&	3.1e{\tt-}04	&	2.0e{\tt-}04 &	1.3e{\tt-}04	&	2.4e{\tt-}04	&	1.9e{\tt-}04	&	1.5e{\tt-}04\\
10935	&	8.8e{\tt-}04	&	NaN	&	7.9e{\tt-}05	&	1.7e{\tt-}04	&	1.1e{\tt-}04           &	8.3e{\tt-}05	&	1.3e{\tt-}04	&	9.9e{\tt-}05	&	8.4e{\tt-}05\\
23436	&	5.1e{\tt-}04	&	4.9e{\tt-}05	&	4.9e{\tt-}05	&	9.2e{\tt-}05	&	6.7e{\tt-}05 &	5.0e{\tt-}05  &	7.5e{\tt-}05	&	6.1e{\tt-}05	&	5.0e{\tt-}05\\
46251	&	3.3e{\tt-}04	&	NaN	&	2.7e{\tt-}05	&	5.3e{\tt-}05	&	4.1e{\tt-}05           &	3.2e{\tt-}05 	&	4.5e{\tt-}05	&	3.6e{\tt-}05	&	3.1e{\tt-}05\\
89372	&	2.1e{\tt-}04	&	NaN	&	1.8e{\tt-}05	&	3.2e{\tt-}05	&	2.5e{\tt-}05           &	2.0e{\tt-}05 	&	2.7e{\tt-}05	&	2.3e{\tt-}05	&	2.0e{\tt-}05\\
\hline
density	&	14.6	&	7.4	&	7.4	&	7.7	&	7.6	& 8.3 &	7.8	&	7.0	&	7.8\\
\hline
\end{tabular}
\caption{Test Problem~\ref{pro11}: RRMS errors $E_{\rm ref}$ for uniform interior grids and projected boundary nodes. 
For meshless finite difference methods we use classical 7-node grid stencil whenever possible.
}
\label{pr11grid}
\end{table}

In the last experiment for Test Problem~\ref{pro11}, we use Halton interior nodes, with boundary nodes obtained by orthogonal
projection. The results can be found in Table~\ref{pr11halton} and Figure~\ref{Pr1_tables}(c), where for \texttt{fem1}
a 3D Delaunay triangulation of $\Xi$ is again created by MATLAB command \texttt{delaunayTriangulation}. 
We use  $\delta=0.9$, $s=1$, $n=3$ and  $k=17$ for \texttt{oct-dist}. In this experiment selection methods \texttt{tet} and \texttt{pQR3}
do not perform well, therefore we do not include them in Figure~\ref{Pr1_tables}(c). 
The extremely large errors of \texttt{pQR3}
for the two finest sets of nodes indicate instability of the system matrix, which will be investigated in detail
for Test Problem~\ref{pro19}. The method \texttt{pQR4} shows significantly
higher errors than the remaining selection algorithms, including the simple selection of 20 nearest points \texttt{20near},
which even looks like the best option in this case. We notice that the mFD method performs very well on
Halton nodes, sometimes even better than on the optimized triangulations, whereas the errors of \texttt{fem1} 
are the highest for this type of $\Xi$.

\begin{table}[htbp!]\scriptsize%
\centering
\renewcommand{\arraystretch}{1.25}
\begin{tabular}{|c||c||c|c|c|c|c||c|c|}
\hline
\multirow{2}{.7cm}{$\# {\Xi _{\rm int}}$} & FEM   
 &	\multicolumn{5}{|c||}{Polyharmonic $r^5$ with quadratic polynomial} & \multicolumn{2}{c|}{Pivoted QR}\\
\cline{2-9}
		& {\scriptsize fem1} & {\scriptsize tet} 	&	{\scriptsize oct} & {\scriptsize oct-dist} 
		& {\scriptsize 20near} & {\scriptsize pQR4sel} & {\scriptsize pQR3} & {\scriptsize pQR4}\\
\hline
303	&	1.7e{\tt-}02	&	NaN	&	3.2e{\tt-}03	&	1.8e{\tt-}03	&	2.0e{\tt-}03	&	1.8e{\tt-}03	&	3.9e{\tt-}03	&	1.0e{\tt-}02\\
627	&	1.1e{\tt-}02	&	NaN	&	1.0e{\tt-}03	&	1.2e{\tt-}03	&	8.3e{\tt-}04	&	1.2e{\tt-}03	&	1.6e{\tt-}03	&	3.0e{\tt-}03\\
1308	&	6.1e{\tt-}03	&	NaN	&	7.8e{\tt-}04	&	9.8e{\tt-}04	&	3.3e{\tt-}04	&	5.0e{\tt-}04	&	9.4e{\tt-}04	&	1.3e{\tt-}03\\
2699	&	3.5e{\tt-}03	&	NaN	&	3.6e{\tt-}04	&	3.3e{\tt-}04	&	2.1e{\tt-}04	&	2.6e{\tt-}04	&	3.8e{\tt-}04	&	7.9e{\tt-}04\\
5493	&	2.1e{\tt-}03	&	NaN	&	2.9e{\tt-}04	&	2.1e{\tt-}04	&	1.4e{\tt-}04	&	1.6e{\tt-}04	&	2.1e{\tt-}04	&	5.1e{\tt-}04\\
11140	&	1.2e{\tt-}03	&	NaN	&	1.9e{\tt-}04	&	1.6e{\tt-}04	&	5.6e{\tt-}05	&	7.8e{\tt-}05	&	1.5e{\tt-}04	&	2.9e{\tt-}04\\
22561	&	7.6e{\tt-}04	&	NaN	&	7.6e{\tt-}05	&	7.1e{\tt-}05	&	2.9e{\tt-}05	&	6.9e{\tt-}05	&	1.2e{\tt-}04	&	1.9e{\tt-}04\\
45513	&	4.3e{\tt-}04	&	NaN	&	3.8e{\tt-}05	&	5.6e{\tt-}05	&	3.5e{\tt-}05	&	4.1e{\tt-}05	&	3.0e{\tt+}08	&	1.1e{\tt-}04\\
91655	&	2.7e{\tt-}04	&	NaN	&	3.0e{\tt-}05	&	2.6e{\tt-}05	&	1.7e{\tt-}05	&	2.1e{\tt-}05	&	2.1e{\tt+}08	&	7.4e{\tt-}05\\
\hline
density	&	16.0	&	16.0	&	16.5	&	16.7	&	19.4	&	19.4	&	9.8	&	19.4 \\
\hline
\end{tabular}
\caption{Test Problem~\ref{pro11}: RRMS error $E_{\rm ref}$  for  
Halton quasi-random interior nodes and projected boundary nodes.
}
\label{pr11halton}
\end{table}

It is well known that the performance of the finite element method depends on the mesh quality, in particular
the shape regularity of the elements, see e.g.\ \cite{FreyGeorge2008}. 
We therefore provide in Table~\ref{pr11meshq} information about the distribution of the aspect ratios
of tetrahedra for each of the three types of triangulations in the above experiments. The
(inverse) \emph{aspect ratio} of a tetrahedron $T$ is given by 
\begin{equation}\label{gamma}
\gamma_T=2\sqrt{6}\,\rho_T/h_T
\end{equation}
where $h_T$ is the diameter and $\rho_T$ is the inradius of the tetrahedron $T$ \cite[p.~1317]{Gmsh_paper09}.
The coefficient $2\sqrt{6}$ is chosen such that $0\le\gamma_T\le1$, with $\gamma_T=0$ for a degenerate
tetrahedron and  $\gamma_T=1$ for a regular tetrahedron. The quality of $T$ is higher when $\gamma_T$ is larger.
Table~\ref{pr11meshq} shows the minimum and the average values of $\gamma_T$ over all tetrahedra of the finest triangulation in
each Table~\ref{pr11fem}, \ref{pr11grid} and \ref{pr11halton}, as well as percentages of tetrahedra with $\gamma_T$
in the ranges 0--0.25, 0.25--0.5, 0.5--0.75 and 0.75--1.0. We compute these statistics using Gmsh \cite{Gmsh} 
for the finest triangulation in each experiments.
We see that the quality of the optimized triangulation used in the first experiment is significantly higher than that 
of the other two triangulations, which explains the smaller errors of \texttt{fem1} on the optimized triangulation. 
In contrast to this, higher shape regularity of the optimized triangulation does not seem to be advantageous for the 
mFD methods.

\begin{table}[htbp!]\scriptsize%
\centering
\renewcommand{\arraystretch}{1.25}
\begin{tabular}{|c||c|c|c|c|c|c|}
\hline
Triangulation	&	$\min\gamma$&	$\avg\gamma$	&	$0<\gamma\leqslant0.25$	&	$0.25<\gamma\leqslant0.5$	
              &	$0.5<\gamma\leqslant0.75$	&	$0.75<\gamma\leqslant1.0$	\\
\hline
Optimized	&	0.50& 0.87	&	0.0\%	&	0.0\% &	6.3\%	&	93.7\%\\
Delaunay/Grid	&	4.6e{\tt-}02& 0.68	&	0.1\%	&	2.2\% &	96.5\%	&	1.2\%	 \\
Delaunay/Halton	&	1.8e{\tt-}03& 0.63	&	4.8\%	&	21.0\% &	45.8\%	&	28.4\%	 \\
\hline
\end{tabular}
\caption{Test Problem~\ref{pro11}: Statistics of the aspect ratio $\gamma_T$ for the simplices  of three 3D 
triangulations used in the tests with finite element method \texttt{fem1}, namely optimized triangulation as in 
Table~\ref{pr11fem},
Delaunay triangulation of gridded nodes  as in Table~\ref{pr11grid}, and Delaunay triangulation of 
Halton nodes  as in Table~\ref{pr11halton}.
}
\label{pr11meshq}
\end{table}

\begin{figure}[!ht]
\begin{center}
\includegraphics[width=7.5cm,height=5.0cm]{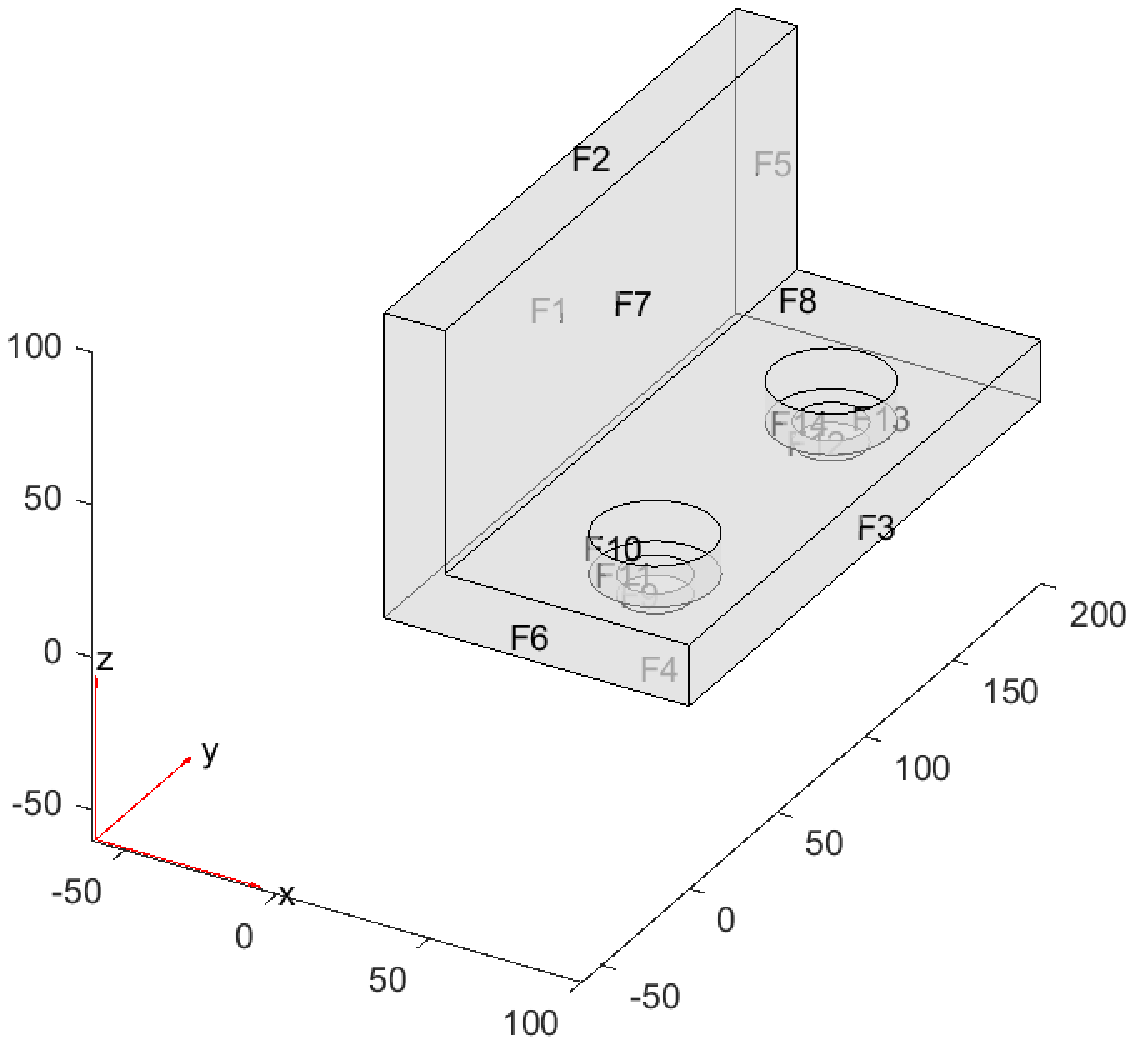}
\caption{Test Problem~\ref{pro19}: Domain \texttt{BracketTwoHoles}; the 3D plot is produced by MATLAB PDE
 Toolbox command {\tt pdegplot}.
}
\label{fig19_1}
\end{center}
\end{figure}

\begin{testproblem}[\texttt{BracketTwoHoles}]
\label{pro19}
Poisson equation $\Delta u = -10$ with zero Dirichlet boundary conditions on the domain $\Omega$ defined in
the STL file  `BracketTwoHoles.stl' shipped with MATLAB  PDE Toolbox \cite{PDEtool},
see Figure~\ref{fig19_1}.
\end{testproblem}

Since the exact solution of this problem is not known, we compute a reference solution by  
the second order finite element method on a triangulation with 1605099 nodes obtained by \texttt{generateMesh} with
$\texttt{Hmax}=1.7$, $\texttt{Hmin}=\texttt{Hmax}/3$,
\texttt{Hgrad}=1.5 and \texttt{GeometricOrder}=2. 

Similar to Test Problem~\ref{pro11}, for our first experiment we generate $\Xi$ as  vertices of optimized triangulations 
produced  by 
\texttt{generateMesh}, with $H_0=8.5$. We use  parameter values $\delta=0.9$, $s=1$, $n=3$ and  $k=13$ for \texttt{oct-dist}.  
The results are presented in Table~\ref{pr19fem}. The errors of the second order finite element method \texttt{fem2} for
optimized triangulations with $H_0=17$ are shown in Table~\ref{pr19fem2}. 
Figure~\ref{Pr2_tables}(a) presents the same errors as functions of nominal  \texttt{nnz}. We see here that 
the versions \texttt{oct-dist}, \texttt{pQR3} and \texttt{pQR4} of mFD produce more accurate solutions than
\texttt{fem1} and compete with the second order method \texttt{fem2}.

\begin{table}[htbp!]\scriptsize %
\centering
\renewcommand{\arraystretch}{1.25}
\begin{tabular}{|c||c||c|c|c|c|c||c|c|}
\hline
\multirow{2}{.7cm}{$\# {\Xi _{\rm int}}$} & FEM   
 &	\multicolumn{5}{|c||}{Polyharmonic $r^5$ with quadratic polynomial} & \multicolumn{2}{c|}{Pivoted QR}\\
\cline{2-9}
		& {\scriptsize fem1}	& {\scriptsize tet} &	{\scriptsize oct} & {\scriptsize oct-dist} 
		& {\scriptsize 20near} & {\scriptsize pQR4sel} & {\scriptsize pQR3} & {\scriptsize pQR4}\\
\hline
871	&	5.5e{\tt-}02	&	4.3e{\tt-}02	&	6.1e{\tt-}02	&	3.2e{\tt-}02	&	7.6e{\tt-}02	&	7.7e{\tt-}02	&	3.8e{\tt-}02	&	2.5e{\tt-}02\\
2086	&	4.3e{\tt-}02	&	NaN	&	5.6e{\tt-}02	&	3.5e{\tt-}02	&	6.3e{\tt-}02	&	6.2e{\tt-}02	&	3.2e{\tt-}02	&	2.2e{\tt-}02\\
4150	&	3.4e{\tt-}02	&	2.9e{\tt-}02	&	3.7e{\tt-}02	&	2.4e{\tt-}02	&	4.4e{\tt-}02	&	4.5e{\tt-}02	&	2.4e{\tt-}02	&	1.6e{\tt-}02\\
9628	&	1.9e{\tt-}02	&	1.8e{\tt-}02	&	2.4e{\tt-}02	&	1.4e{\tt-}02	&	3.3e{\tt-}02	&	3.2e{\tt-}02	&	1.2e{\tt-}02	&	1.1e{\tt-}02\\
19493	&	1.5e{\tt-}02	&	1.3e{\tt-}02	&	2.0e{\tt-}02	&	1.2e{\tt-}02	&	2.4e{\tt-}02	&	2.2e{\tt-}02	&	1.1e{\tt-}02	&	9.5e{\tt-}03\\
41215	&	8.6e{\tt-}03	&	NaN	&	1.1e{\tt-}02	&	6.2e{\tt-}03	&	1.5e{\tt-}02	&	1.5e{\tt-}02	&	7.8e{\tt-}03	&	5.2e{\tt-}03\\
86699	&	6.8e{\tt-}03	&	NaN	&	9.3e{\tt-}03	&	6.0e{\tt-}03	&	1.1e{\tt-}02	&	1.1e{\tt-}02	&	6.1e{\tt-}03	&	4.4e{\tt-}03\\
178432	&	4.5e{\tt-}03	&	NaN	&	7.2e{\tt-}03	&	4.7e{\tt-}03	&	6.8e{\tt-}03	&	7.5e{\tt-}03	&	4.9e{\tt-}03	&	3.1e{\tt-}03\\
\hline
density	&	14.2	&	14.2	&	15.9	&	12.4	&	18.8	&	18.7	&	9.5	&	18.7 \\
\hline
\end{tabular}
\caption{Test Problem~\ref{pro19}:  RRMS errors $E_{\rm ref}$ for optimized
triangulations with $H_0=8.5$.
}
\label{pr19fem}
\end{table}

\begin{table}[htbp!]\scriptsize %
\centering
\renewcommand{\arraystretch}{1.25}
\begin{tabular}{|c|c|c|c|c|c|c|c|c||c|}
\hline
$\# {\Xi _{\rm int}}$: & 1197& 2647& 	4121& 	9507& 	21034& 	39659& 	87900& 	173262 & density\\
\hline
fem2:	&	4.0e{\tt-}02	&	2.3e{\tt-}02	&	1.6e{\tt-}02	&	9.7e{\tt-}03	&	6.6e{\tt-}03	&	5.0e{\tt-}03	&	3.1e{\tt-}03	&	2.3e{\tt-}03	&	25.9\\
\hline
\end{tabular}
\caption{Test Problem~\ref{pro19}:  RRMS errors $E_{\rm ref}$ of the quadratic finite element method 
for optimized triangulations with $H_0=17$.}
\label{pr19fem2}
\end{table}

In addition to solving mFD linear system \eqref{prd} by using MATLAB's backslash operator,
we run an iterative solver for all successful selection methods of Table~\ref{pr19fem}. Following \cite{BFFB17}
we use the BiCGSTAB iterative method  with default tolerance $10^{-6}$ and 
maximum number of iterations \texttt{maxit}=1000, applying the incomplete LU factorization ILU(0) as
preconditioner and reverse Cuthill-McKee ordering. The errors of the numerical solutions obtained this way are the same as in Table~\ref{pr19fem}. The number of
iterations ($\#$iter) before  BiCGSTAB terminates is shown in Table~\ref{pr19fem_iter} for each selection method except
\texttt{tet} and each node set.
Note that the number of iterations is sometimes a half-integer because each BiCGSTAB iteration consists of two steps: a
biconjugate  gradient step alternating with a GMRES step for additional stability. We observe that $\#$iter remains small for
all methods even on the finest set of 178432 nodes, which indicates that the simplest and cheapest preconditioning method 
ILU(0) works very well for the mFD system matrices.

\begin{table}[htbp!]\scriptsize %
\centering
\renewcommand{\arraystretch}{1.25}
\begin{tabular}{|c||c|c|c|c|c|c|}
\hline
\multirow{2}{.7cm}{$\# {\Xi _{\rm int}}$}   
 &	\multicolumn{6}{|c|}{Number of iterations $\#$iter}\\
\cline{2-7}
		& 	{\scriptsize oct} & {\scriptsize oct-dist} 
		& {\scriptsize 20near} & {\scriptsize pQR4sel} & {\scriptsize pQR3} & {\scriptsize pQR4}\\
\hline
871	&	4	&	3.5	&	3	&	3.5	&	3.5	&	3\\
2086	&	5	&	4	&	4	&	4	&	4.5	&	3.5\\
4150	&	5.5	&	5	&	4.5	&	5	&	5.5	&	4.5\\
9628	&	7	&	7	&	7.5	&	6.5	&	7.5	&	5.5\\
19493	&	10	&	9	&	8	&	8.5	&	9.5	&	7.5\\
41215	&	12.5	&	11.5	&	10.5	&	11.5	&	13.5	&	10\\
86699	&	14.5	&	13.5	&	13.5	&	12.5	&	14.5	&	14\\
178432	&	19.5	&	21	&	16.5	&	19	&	18.5	&	16\\
\hline
\end{tabular}
\caption{Test Problem~\ref{pro19}:  Number of iterations of BiCGSTAB for optimized
triangulations as in Table~\ref{pr19fem}.
}
\label{pr19fem_iter}
\end{table}

In the next experiment we generate $\Xi$ as vertices of an unoptimized triangulation of $\Omega$. Since this domain
is  not convex, \texttt{delaunayTriangulation} is not available to triangulate it, and we use Gmsh with 
smoothing and optimization switched off, as explained above, with $H_0=8.5$. The quality of both triangulations is compared 
in Table~\ref{pr19meshq}. Parameters of \texttt{oct-dist} are chosen as $\delta=0.9$, $s=3$, $n=6$ and  $k=17$. 
The results are presented in Table~\ref{pr19gmsh} and Figure~\ref{Pr2_tables}(b). The only selection methods included in
Figure~\ref{Pr2_tables}(b) are \texttt{oct-dist}, \texttt{pQR4} and \texttt{pQR4sel} because \texttt{tet}, \texttt{oct} and
\texttt{20near} fail to generate sets of influence with weights satisfying the exactness condition \eqref{exap}, whereas
\texttt{pQR3} too often produces high errors and even leads to a singular system matrix of \eqref{prd} on the finest node set.
Nevertheless, the error of \texttt{pQR3} for the discretization  with $\#\Xi_{\rm int}=86541$ is among the best, with by far 
lower density of the system matrix than in any other method in the table.
In fact, the error is also poor for \texttt{pQR4} on one of the sets $\Xi$, see the peak in its plot 
in Figure~\ref{Pr2_tables}(b), but we
nevertheless included \texttt{pQR4} in the figure because its errors for the other node sets are the best compared to other methods.

\begin{table}[htbp!]\scriptsize%
\centering
\renewcommand{\arraystretch}{1.25}
\begin{tabular}{|c||c|c|c|c|c|c|}
\hline
Triangulation	&	$\min\gamma$&	$\avg\gamma$	&	$0<\gamma\leqslant0.25$	&	$0.25<\gamma\leqslant0.5$	
              &	$0.5<\gamma\leqslant0.75$	&	$0.75<\gamma\leqslant1.0$	\\
\hline
Optimized	&	0.51& 0.88	&	0.0\%  	&	0.0\% &	2.8\%	&	97.2\%\\
Unoptimized	&	3.1e{\tt-}05& 0.78	&	1.7\%	&	5.8\% &	25.1\%	&	67.4\%	 \\
\hline
\end{tabular}
\caption{Test Problem~\ref{pro19}: Statistics of the aspect ratio $\gamma_T$ for the simplices  of two types of 3D 
triangulations used in the tests with finite element method \texttt{fem1}, namely an optimized triangulation as in 
Table~\ref{pr19fem}, and an unoptimized one as in Table~\ref{pr19gmsh}.
}
\label{pr19meshq}
\end{table}

\begin{table}[htbp!]\scriptsize %
\centering
\renewcommand{\arraystretch}{1.25}
\begin{tabular}{|c||c||c|c|c|c|c|c|c||c|c|}
\hline
\multirow{2}{.7cm}{$\# {\Xi _{\rm int}}$} & FEM   
 &	\multicolumn{7}{|c||}{Polyharmonic $r^5$ with quadratic polynomial} & \multicolumn{2}{c|}{Pivoted QR}\\
\cline{2-11}
		& {\scriptsize fem1}	& {\scriptsize tet} &	{\scriptsize oct} & {\scriptsize oct-dist} 
		& {\scriptsize 20near} & {\scriptsize 30near}& {\scriptsize 40near}& {\scriptsize pQR4sel} & {\scriptsize pQR3} & {\scriptsize pQR4}\\
\hline
726	&	1.3e{\tt-}01	&	NaN	&	NaN	&	7.1e{\tt-}02	&	NaN	&	NaN	&	NaN	&	1.1e{\tt-}01	&	2.1e{\tt-}01	&	1.1e{\tt-}01	\\
1491	&	1.1e{\tt-}01	&	NaN	&	NaN	&	7.1e{\tt-}02	&	NaN	&	NaN	&	1.2e{\tt-}01	&	8.3e{\tt-}02	&	1.8e{\tt+}00	&	4.5e{\tt-}02\\
2529	&	6.8e{\tt-}02	&	NaN	&	NaN	&	4.9e{\tt-}02	&	NaN	&	NaN	&	NaN	&	5.8e{\tt-}02	&	1.2e{\tt-}01	&	3.4e{\tt-}02	\\
4954	&	5.0e{\tt-}02	&	NaN	&	NaN	&	3.6e{\tt-}02	&	NaN	&	NaN	&	NaN	&	4.1e{\tt-}02	&	3.6e{\tt-}02	&	1.6e{\tt-}01	\\
9940	&	3.4e{\tt-}02	&	NaN	&	NaN	&	2.6e{\tt-}02	&	NaN	&	NaN	&	5.4e{\tt-}02	&	3.2e{\tt-}02	&	7.3e{\tt-}02	&	1.5e{\tt-}02	\\
19834	&	2.4e{\tt-}02	&	NaN	&	NaN	&	1.7e{\tt-}02	&	NaN	&	NaN	&	NaN	&	2.2e{\tt-}02	&	6.7e{\tt-}02	&	1.3e{\tt-}02	\\
41941	&	1.6e{\tt-}02	&	NaN	&	NaN	&	1.3e{\tt-}02	&	NaN	&	NaN	&	3.0e{\tt-}02	&	1.5e{\tt-}02	&	9.9e{\tt-}03	&	8.1e{\tt-}03	\\
86541	&	1.1e{\tt-}02	&	NaN	&	1.0e{\tt-}02	&	9.6e{\tt-}03	&	NaN	&	NaN	&	2.1e{\tt-}02	&	1.1e{\tt-}02	&	6.7e{\tt-}03	&	5.9e{\tt-}03	\\
178450	&	7.4e{\tt-}03	&	NaN	&	NaN	&	6.7e{\tt-}03	&	NaN	&	NaN	&	NaN	&	7.7e{\tt-}03	&	Inf	&	3.9e{\tt-}03	\\
\hline
density	&	14.9	&	14.9	&	15.9	&	15.9	&	18.4	&	27.4	&	36.3	&	18.8	&	9.6	&		18.8\\
\hline
\end{tabular}
\caption{Test Problem~\ref{pro19}: RRMS errors $E_{\rm ref}$ for unoptimized triangulations.
}
\label{pr19gmsh}
\end{table}

Since the pQR method produces accurate numerical differentiation formulas \cite{D19arxiv}, the reasons for the large 
errors seen on certain node sets seem to be related to the global properties of the linear system \eqref{prd}. 
In order to investigate this, we compute the \emph{stability constant}
\begin{equation}\label{stab}
\sigma:=\texttt{condest}(A^T)/\|A\|_\infty,
\end{equation}
where $A$ is the system matrix of \eqref{prd}. This formula, suggested in 
\cite{Schaback_error_analysis17}, uses an efficient estimator of the 1-norm condition number available
in MATLAB via command \texttt{condest}. Therefore $\sigma$ is a good approximation of the norm 
$\|A^{-1}\|_\infty$ that measures stability in the classical error analysis of the Finite Difference Method for elliptic
problems. Unfortunately, there are no theoretical results that would allow us to estimate the stability constant, and the error analysis
of the meshless finite difference methods is underdeveloped \cite{D19arxiv2}.

We present the stability constant $\sigma$ and the number of BiCGSTAB iterations $\#$iter for the 
selection methods \texttt{oct-dist, pQR4sel, pQR3} and  \texttt{pQR4} in Table~\ref{pr19gmsh_sigma}. 
We see that large errors are always associated with
relatively high values of  $\sigma$. For the two most stable methods \texttt{oct-dist} and \texttt{pQR4sel}, with monotone graphs in
Figure~\ref{Pr2_tables}(b), we always see $\sigma< 100$, whereas \texttt{pQR4} produces matrices with  $\sigma\ge 300$
for $i\le 5$, and $\sigma>3000$ for the peak of Figure~\ref{Pr2_tables}(b) at $i=3$. The large errors of \texttt{pQR3} 
for $i\le 5$ and $i=8$ correspond to extremely large $\sigma$, and very good results for  $i\in \{6,7\}$ correspond to 
$\sigma\approx 61$, close to the minimum of $\sigma=58.6$ in Table~\ref{pr19gmsh_sigma}. 
Note that $56.3\le\sigma\le 68.9$ in all mFD tests for optimized triangulations reported in Table~\ref{pr19fem}.
Iteration counts for \texttt{oct-dist} and \texttt{pQR4sel} for all node sets are close to those in Table~\ref{pr19fem_iter} 
for the optimized triangulation. The same is true for \texttt{pQR3} when $i\in \{6,7\}$ and \texttt{pQR4} when 
$i\in \{6,7,8\}$, with small $\sigma$ and very good errors in Table~\ref{pr19gmsh}. Higher iteration numbers for \texttt{pQR4}
when $0\le i\le 4$ are associated with large $\sigma$, and the failure of the iteration (indicated as `fail' in the table) only
occurs when the stability constant is extremely large, $\sigma>10^{15}$.

\begin{table}[htbp!]\scriptsize%
\centering
\renewcommand{\arraystretch}{1.25}
\begin{tabular}{|c||c|c||c|c||c|c||c|c|}
\hline
\multirow{2}{.7cm}{$\# {\Xi _{\rm int}}$} & \multicolumn{2}{|c||}{oct-dist}  & \multicolumn{2}{|c||}{pQR4sel} 
& \multicolumn{2}{|c||}{pQR3} & \multicolumn{2}{|c|}{pQR4}\\
\cline{2-9} 
& $\sigma$ & $\#$iter & $\sigma$ & $\#$iter & $\sigma$ & $\#$iter & $\sigma$ & $\#$iter \\
\hline
 726	&	91.4	&	6.5	&	76.1	&	7	&	2.9e{\tt +}16	&	fail	&	1.6e{\tt +}03	&	36.5\\
1491	&	58.6	&	24	&	82.9	&	37.5	&	1.6e{\tt +}17	&	fail	&	3.4e{\tt +}03	&	477\\
2529	&	59.0	&	11.5	&	68.2	&	19.5	&	1.6e{\tt +}17	&	fail	&	585.8	&	157\\
4954	&	58.7	&	15.5	&	64.9	&	23.5	&	2.1e{\tt +}17	&	fail	&	3.5e{\tt +}03	&	46.5\\
9940	&	59.2	&	31	&	58.6	&	10	&	2.5e{\tt +}16	&	fail	&	703.3	&	42.5\\
19834	&	59.8	&	9.5	&	59.6	&	10	&	6.6e{\tt +}15	&	fail	&	304.2	&	12\\
41941	&	60.2	&	10.5	&	60.3	&	12	&	60.6	&	13.5	&	60.7	&	10\\
86541	&	60.6	&	13.5	&	60.5	&	13.5	&	60.8	&	16	&	60.7	&	13\\
178450	&	60.7	&	17	&	60.7	&	18	&	Inf	&	fail	&	61.1	&	15.5\\
\hline
\end{tabular}
\caption{Test Problem~\ref{pro19}: stability constant $\sigma$ and number of BiCGSTAB  iterations $\#$iter for 
unoptimized triangulations as in Table~\ref{pr19gmsh}. We write 'fail' for $\#$iter if BiCGSTAB  terminates
unsuccessfully.
}
\label{pr19gmsh_sigma}
\end{table}

\begin{figure}[!ht]
\begin{center}
\hbox{\hspace*{-40pt}\subfigure[STL triangulation of the boundary]{\includegraphics[width=8.5cm]{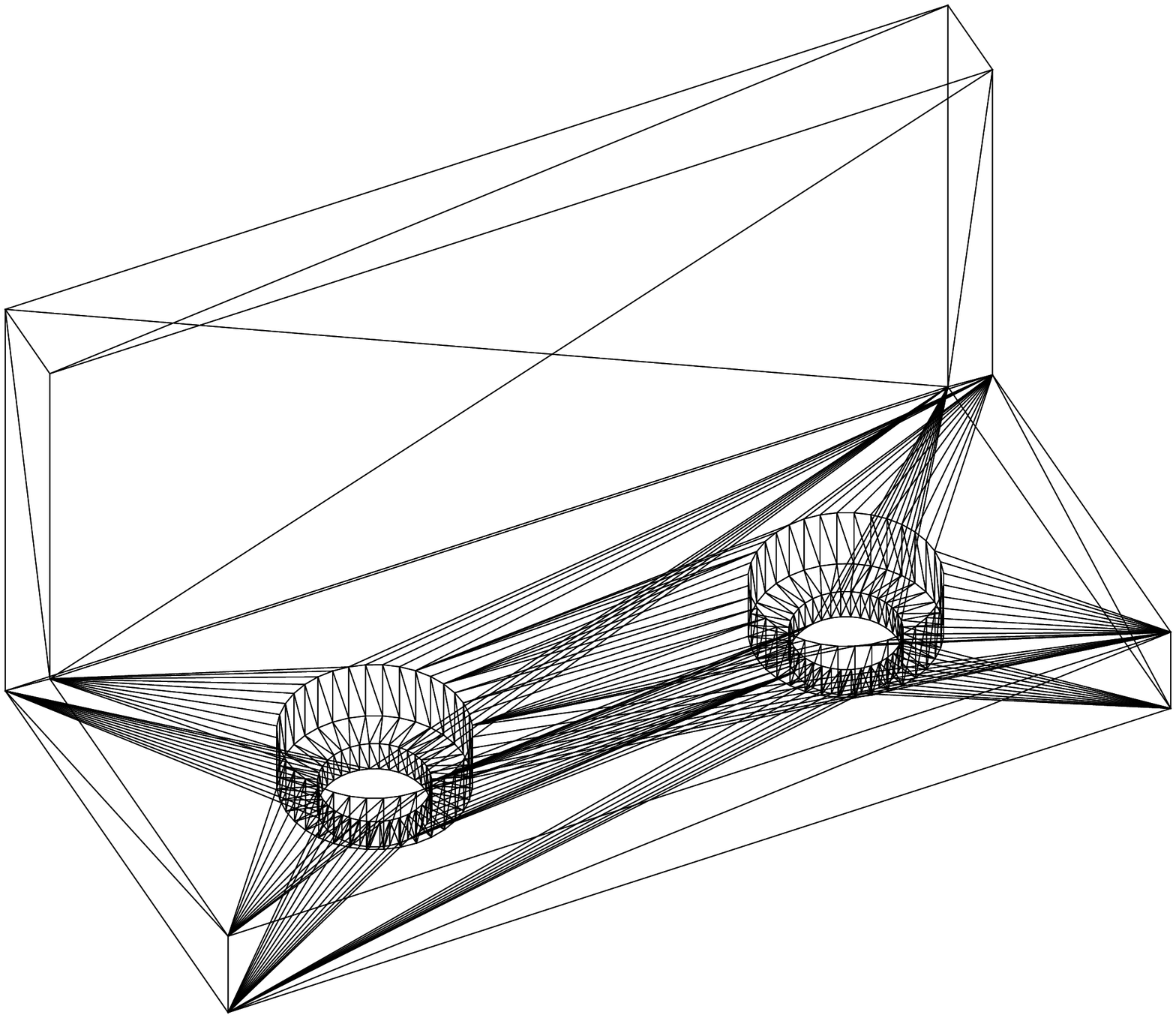}}
\hspace*{-40pt}\subfigure[30 neighbors in bad position]{\includegraphics[width=10cm]{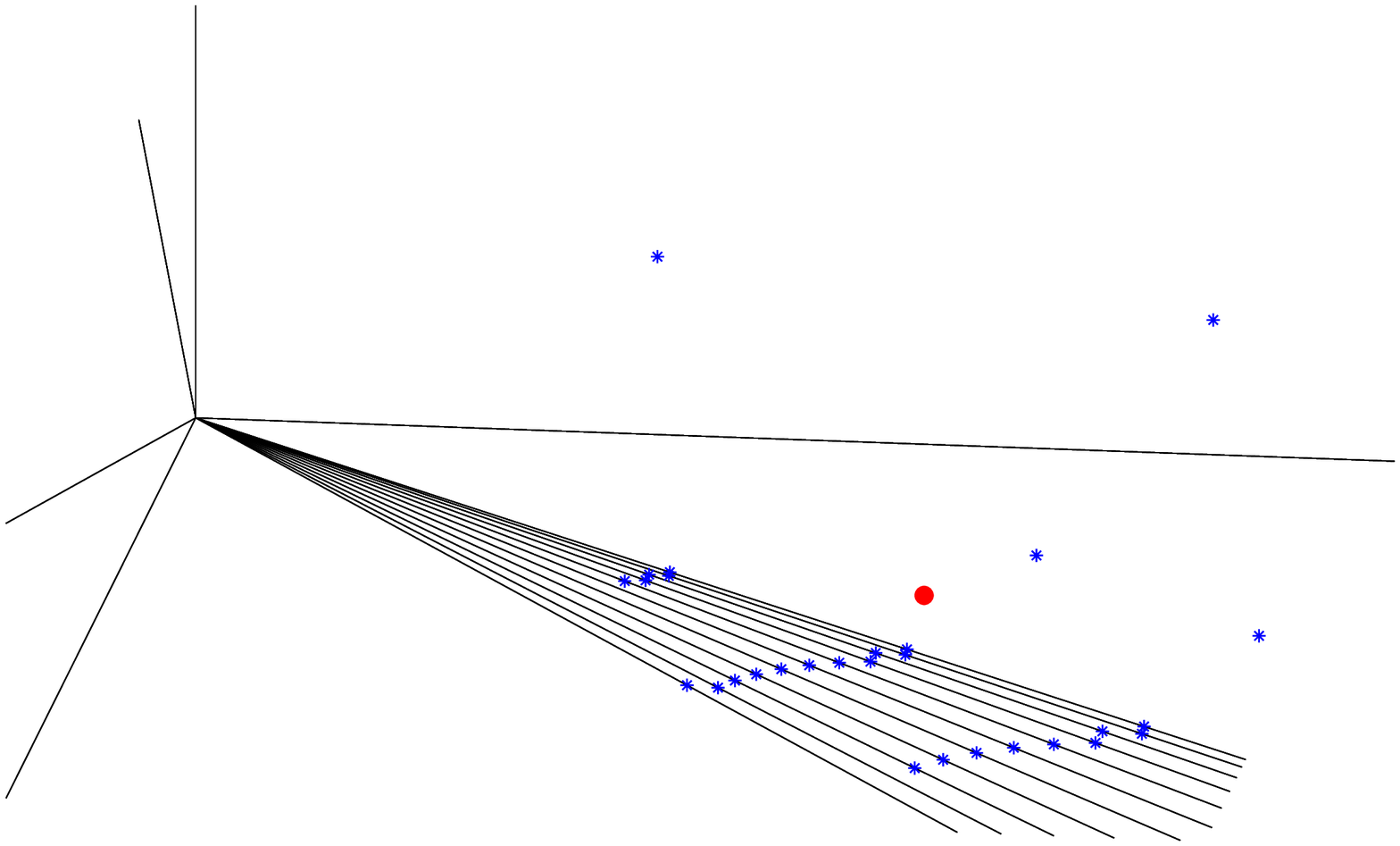}}}
\caption{Test Problem~\ref{pro19}: (a) STL model and (b) 30 nearest vertices of the unoptimized triangulation
for the node marked as a red dot. Note that all nodes shown are lying in the horizontal plane except of the red node itself and 
the two topmost nodes.
\label{Pr2_stl}}
\end{center}
\end{figure}

\begin{figure}[!ht]
\begin{center}
\hbox{\subfigure[Optimized triangulation]{\includegraphics[width=7.5cm,height=5.0cm]{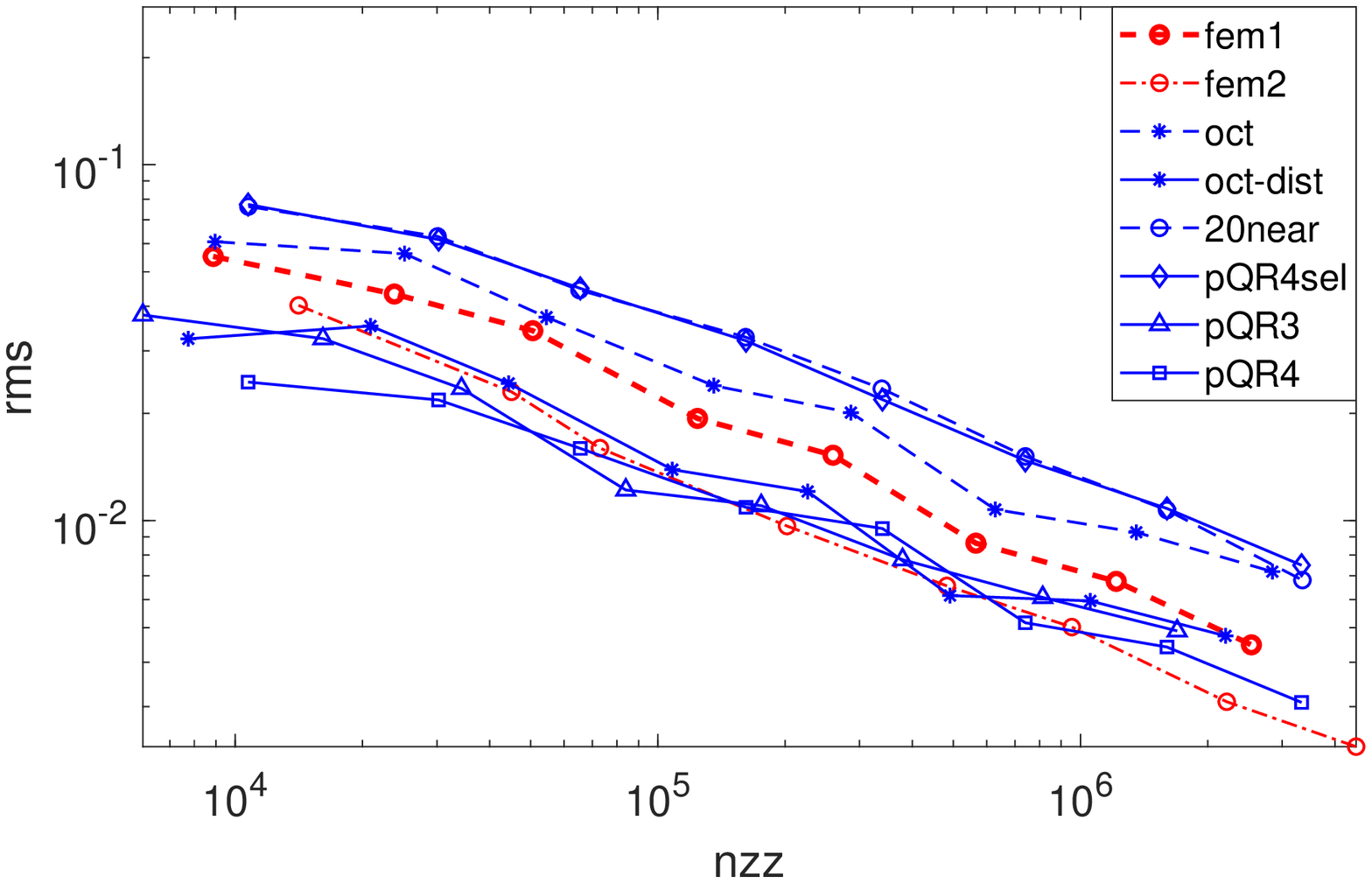}}
\subfigure[Unoptimized triangulation]{\includegraphics[width=7.5cm,height=5.0cm]{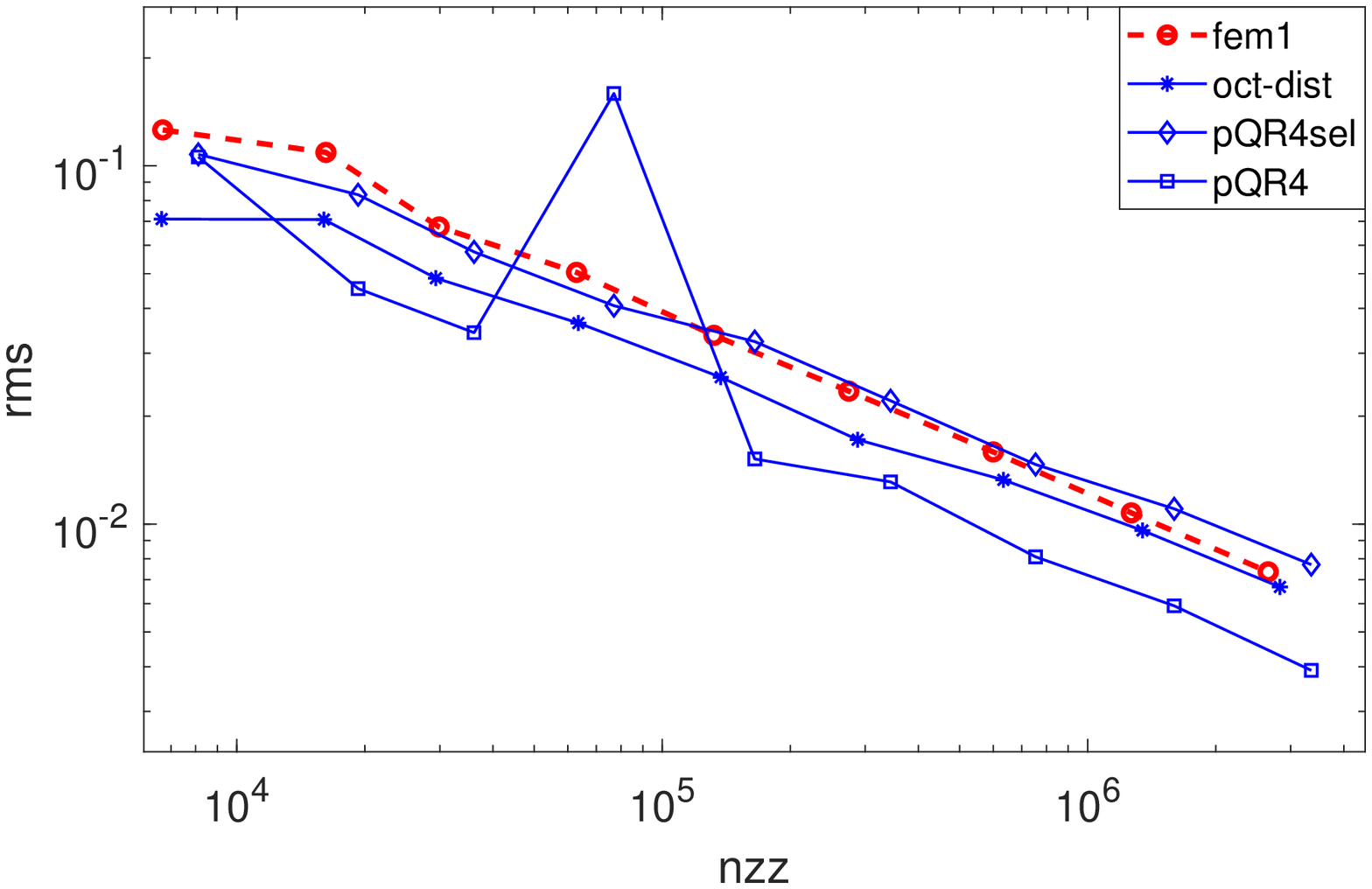}}}
\hbox{\subfigure[Interior grid and boundary nodes of (a)]{\includegraphics[width=7.5cm,height=5.0cm]{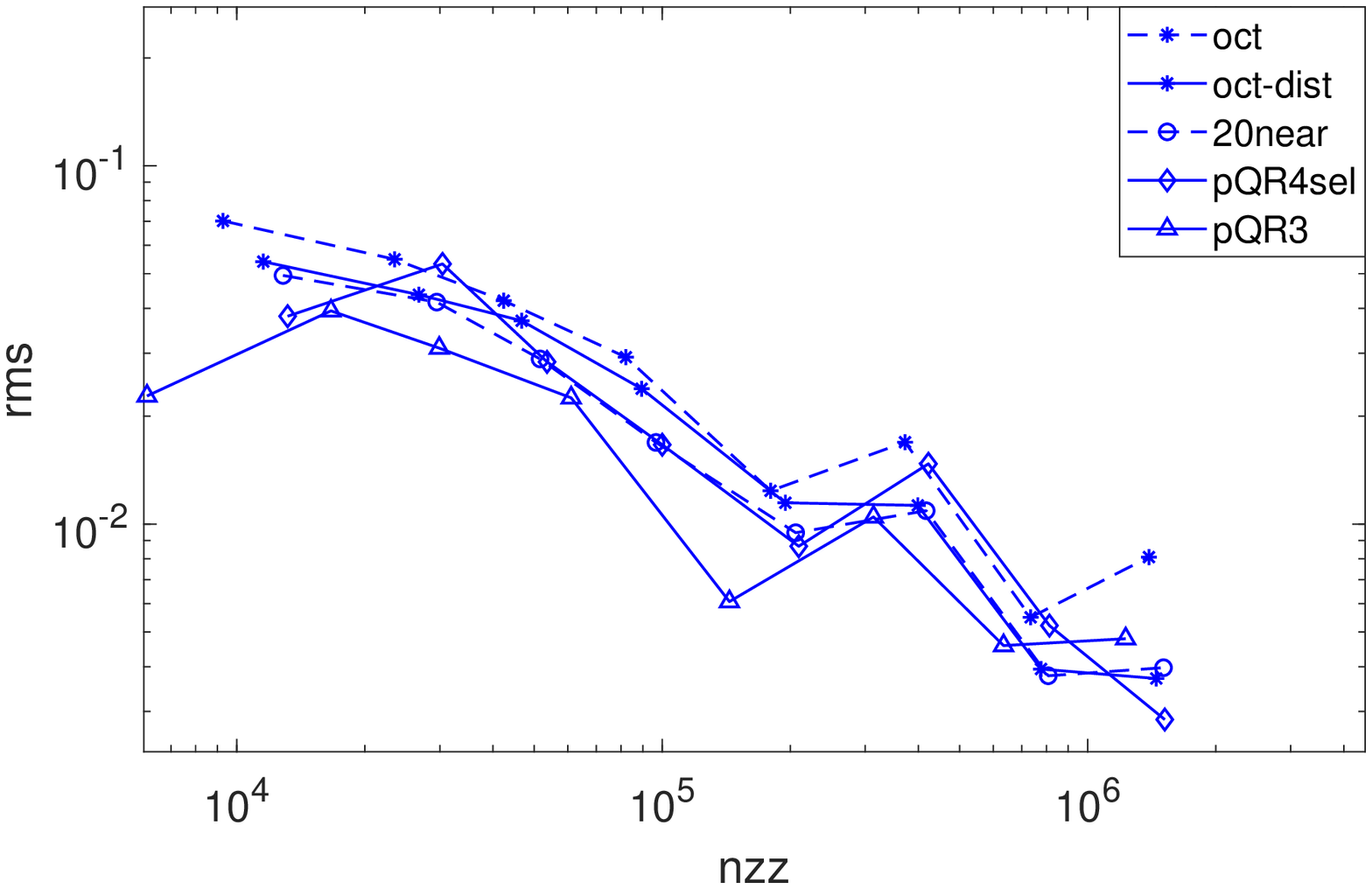}}
\subfigure[Halton interior and projected boundary nodes]{\includegraphics[width=7.5cm,height=5.0cm]{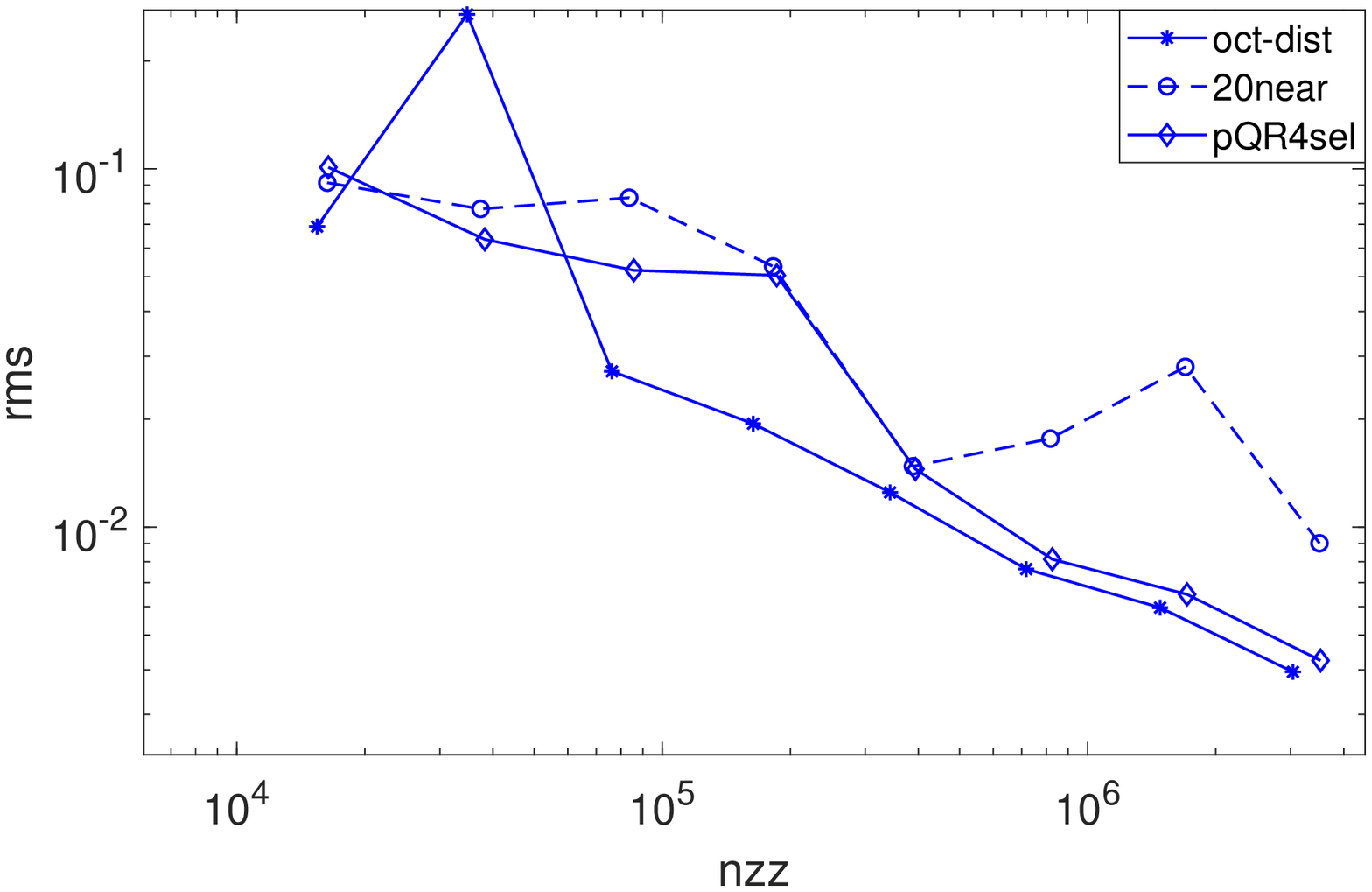}}}
\caption{Test Problem~\ref{pro19}: RRMS errors as functions of nominal \texttt{nnz}, corresponding to 
(a)  Tables~\ref{pr19fem} and \ref{pr19fem2}, (b) Table~\ref{pr19gmsh}, (c) Table~\ref{pr19gridfem}, 
(d) Table~\ref{pr19project}
\label{Pr2_tables}}
\end{center}
\end{figure}

Since \texttt{20near} fails completely, we tried to increase the number of nearest neighbors,
and included in Table~\ref{pr19gmsh} the errors for \texttt{30near} and \texttt{40near} even
if the density of the system matrix gets much higher.
We see that the results remain very poor even for 40 nearest neighbors. This may be related to the fact that STL models
often include many thin and long triangles on the boundary, and  typical triangulation algorithms 
start by discretizing boundary edges, then boundary triangles, and only after that proceed to adding vertices in the 
interior of the domain. As a result, neighborhoods of certain interior nodes include many coplanar boundary nodes 
that are not favorable for numerical differentiation, see Figure~\ref{Pr2_stl}. Therefore it is helpful when the selection 
algorithm abandons some of these coplanar nodes even if they are the closest neighbors to $\zeta$.

In the next experiment we generate nodes by a uniform grid in the interior and orthogonal projection to the boundary,
again using the classical 7-node grid stencil whenever possible.
Unfortunately, no selection method shows a good performance in this case. The methods \texttt{20near}, \texttt{oct} and
\texttt{oct-dist} fail to provide polynomial exactness \eqref{exap}. The errors of \texttt{30near}, \texttt{pQR4sel},
and \texttt{pQR3} are shown in Table~\ref{pr19gridpro} and are disappointing. The performance of
\texttt{pQR4} is not reported, as it is worse than that of  \texttt{pQR4sel} and \texttt{pQR3}.
The data in Table~\ref{pr19gridpro} show that situations where the error of the mFD method is particularly bad correspond to 
stability constant several magnitudes higher than the `normal' values of around 60. On the other hand, whenever $\sigma$ is in
the 60s, the performance of \texttt{pQR4sel} and \texttt{pQR3} is quite competitive in comparison to the nodes from
optimized triangulations. The errors of \texttt{30near} even in these cases are significantly larger. The iterative method
performs well whenever $\sigma$ is small, but often fails already for $\sigma$ approaching 1000.

\begin{table}[htbp!]\scriptsize%
\centering
\renewcommand{\arraystretch}{1.25}
\begin{tabular}{|c||c|c|c||c|c|c||c|c|c|}
\hline
\multirow{2}{.7cm}{\scriptsize $\# {\Xi _{\rm int}}$}	& \multicolumn{3}{|c||}{\scriptsize 30near}
& \multicolumn{3}{|c||}{\scriptsize pQR4sel} & \multicolumn{3}{|c|}{\scriptsize pQR3} 
	\\
\cline{2-10}	
&  $E_{\rm ref}$ & $\sigma$ &$\#$iter & $E_{\rm ref}$ & $\sigma$ &$\#$iter & $E_{\rm ref}$ & $\sigma$ & $\#$iter\\
\hline
1004	&	1.3e{\tt-}01	&	67.4         & 3.5   &1.2e{\tt-}01  &  617.9          &6  & 1.2e{\tt-}01  & 884.1         & fail  \\
2572	&	1.1e{\tt-}01	& 68.1         & 5     &6.0e{\tt-}02  &  63.5           &5.5  & 4.7e{\tt-}02  & 62.1          & 5     \\
4344	&	2.0e{\tt-}01	&	2.1e{\tt+}03 & 437.5 &4.4e{\tt-}01  &  9.0e{\tt+}03   &fail & 2.0e{\tt-}01  & 4.1e{\tt+}03  & fail  \\
8883	&	4.6e{\tt-}01	&	5.9e{\tt+}03 & fail  &5.1e{\tt-}01  &  1.2e{\tt+}04   &fail & 5.9e{\tt-}01  & 5.4e{\tt+}03  & fail  \\
21228	&	2.8e{\tt-}02	&	62.9         & 9.5    &1.6e{\tt-}02  &  61.8           &10   & 1.5e{\tt-}02  & 362.4         & 18  \\
45808	&	3.6e{\tt-}02	&64         & 13    &1.7e{\tt-}02  &  62.3         &12.5   & 1.1e{\tt-}02  & 61.7          & 14  \\
92528	&	2.0e{\tt-}02	&	62.9         & 17    &1.7e{\tt-}02  &  62.8           &16 & 6.7e{\tt-}03  & 61.9          & 16.5  \\
176179&	6.8e{\tt-}01  & 7.8e{\tt+}03 & fail  &7.8e{\tt-}01  &  3.5e{\tt+}04   &fail & 7.3e{\tt-}01  & 1.5e{\tt+}04  & fail  \\
\hline
density	&	\multicolumn{3}{|c||}{10.5}	&	\multicolumn{3}{|c||}{7.1}&	\multicolumn{3}{|c|}{6.9}	\\
\hline
\end{tabular}
\caption{Test Problem~\ref{pro19}: RRMS error $E_{\rm ref}$, stability constant $\sigma$ and the number of BiCGSTAB
iterations  $\#$iter for uniform interior grids and projected boundary nodes. The 7-node grid stencil \eqref{7star} is used
whenever possible.
}
\label{pr19gridpro}
\end{table}

The performance of mFD is much better if we combine interior grid nodes with the discretization of the boundary taken from
the optimized triangulation of the first experiment for this test problem. More precisely, we choose $H_0=8.5$ 
and generate, for each $i=0,\ldots,7$, the interior nodes as a grid with target spacing $h=0.9\,\texttt{Hmax}$, where
$\texttt{Hmax}=H_02^{-i/3}$. To obtain boundary nodes, we run \texttt{generateMesh} with the above $\texttt{Hmax}$, and use the
boundary vertices of the resulting triangulation as boundary nodes of $\Xi$. The results are presented in 
Table~\ref{pr19gridfem} and Figure~\ref{Pr2_tables}(c). Parameters of \texttt{oct-dist} are 
$\delta=0.9$, $s=1$, $n=3$ and $k=18$, as in the case of interior grid nodes for Test Problem~\ref{pro11}. The performance of 
mFD with all selection methods but \texttt{pQR4} is good, although the error graphs in Figure~\ref{Pr2_tables}(c)
are not monotone. In particular, the set $\Xi$ with 45808 interior nodes ($i=5$) seems unfavorable for all methods.
Nevertheless, looking for example at the graph for \texttt{pQR3}, we see that the point corresponding to $i=4$ that makes this
graph looking so irregular, in fact shows an exceptionally good error, by far the best for comparable number of interior
nodes (around 20000) in all experiments for Test Problem~\ref{pro19}. The stability constant $\sigma$ in all tests in
Table~\ref{pr19gridfem} except \texttt{pQR4} satisfies $58.2\le\sigma \le 67.4$. For \texttt{pQR4} we have $\sigma=90.4$ when $i=2$ and 413.0 when $i=5$, otherwise 
$\sigma$ is also below 70, which fits well into the data in Table~\ref{pr19gridfem}, where the errors for $i=2$ and 5 are
significantly higher for \texttt{pQR4} than with other selection methods. In addition to \texttt{20near}, we included in 
Table~\ref{pr19gridfem} the errors of mFD for selection methods \texttt{30near} and \texttt{40near}. However, in contrast to 
Test Problem~\ref{pro11}, this does not improve the error in comparison to \texttt{20near}.
The iteration numbers of BiCGSTAB are presented in Table~\ref{pr19gridfem_iter}. They are slightly higher than those for the
nodes of the optimized triangulations in Table~\ref{pr19fem_iter}, but are still below 26 even on the finest node set.

\begin{table}[htbp!]\scriptsize %
\centering
\renewcommand{\arraystretch}{1.25}
\begin{tabular}{|c||c|c|c|c|c|c||c|c|}
\hline
\multirow{2}{.7cm}{$\# {\Xi _{\rm int}}$} &\multicolumn{6}{c||}{Polyharmonic $r^5$ with quadratic polynomial} & \multicolumn{2}{c|}{Pivoted QR}\\
\cline{2-9}
& {\scriptsize oct} & {\scriptsize oct-dist} & {\scriptsize 20near} & {\scriptsize 30near} & {\scriptsize 40near} & {\scriptsize pQR4sel}   & {\scriptsize pQR3} & {\scriptsize pQR4}\\
\hline
1004	  &	7.0e{\tt-}02	&	5.4e{\tt-}02	&	4.9e{\tt-}02 & 6.7e{\tt-}02	 & 6.4e{\tt-}02 &	3.8e{\tt-}02 	&	2.3e{\tt-}02	&	7.0e{\tt-}02\\
2572	  &	5.5e{\tt-}02	&	4.4e{\tt-}02	&	4.2e{\tt-}02 & 5.5e{\tt-}02	 & 4.9e{\tt-}02 &	5.3e{\tt-}02 	&	3.9e{\tt-}02	&	2.7e{\tt-}01\\
4344	  &	4.2e{\tt-}02	&	3.7e{\tt-}02	&	2.9e{\tt-}02 & 3.8e{\tt-}02	 & 3.7e{\tt-}02 &	2.8e{\tt-}02 	&	3.1e{\tt-}02	&	3.7e{\tt-}02\\
8883  	&	2.9e{\tt-}02	&	2.4e{\tt-}02	&	1.7e{\tt-}02 & 2.5e{\tt-}02	 & 3.0e{\tt-}02 &	1.7e{\tt-}02 	&	2.3e{\tt-}02	&	2.7e{\tt-}02\\
21228  	&	1.2e{\tt-}02	&	1.2e{\tt-}02	&	9.5e{\tt-}03 & 1.2e{\tt-}02	 & 1.3e{\tt-}02	&	8.7e{\tt-}03  &	6.1e{\tt-}03	&	2.0e{\tt-}02\\
45808	  &	1.7e{\tt-}02	&	1.1e{\tt-}02	&	1.1e{\tt-}02 & 1.4e{\tt-}02	 & 1.2e{\tt-}02	&	1.5e{\tt-}02  &	1.1e{\tt-}02	&	7.7e{\tt-}02\\
92528	  &	5.5e{\tt-}03	&	3.9e{\tt-}03	&	3.8e{\tt-}03 & 4.8e{\tt-}03	 & 4.8e{\tt-}03	&	5.2e{\tt-}03  &	4.6e{\tt-}03	&	1.8e{\tt-}02\\
176179	&	8.1e{\tt-}03	&	3.7e{\tt-}03	&	4.0e{\tt-}03&	3.8e{\tt-}03 & 1.2e{\tt-}02	&	2.9e{\tt-}03 	&	4.8e{\tt-}03	&	1.1e{\tt-}02\\
\hline
density	&	7.9	&	8.2	&	8.6	&	9.9	&	11.4	&	8.6	&	7.0	&	8.6\\
\hline
\end{tabular}
\caption{Test Problem~\ref{pro19}: RRMS error  $E_{\rm ref}$ for uniform interior grids and 
boundary nodes from optimized triangulations. The 7-node grid stencil \eqref{7star} is used whenever possible.
}
\label{pr19gridfem}
\end{table}

\begin{table}[htbp!]\scriptsize %
\centering
\renewcommand{\arraystretch}{1.25}
\begin{tabular}{|c||c|c|c|c|c|c|c|c|}
\hline
\multirow{2}{.7cm}{$\# {\Xi _{\rm int}}$}   
 &	\multicolumn{8}{|c|}{Number of iterations}\\
\cline{2-9}
		& 	{\scriptsize oct} & {\scriptsize oct-dist} 
		& {\scriptsize 20near}& {\scriptsize 30near}& {\scriptsize 40near} & {\scriptsize pQR4sel} & {\scriptsize pQR3} & {\scriptsize pQR4}\\
\hline
1004	&	4.5	&	3.5	&	3.5	&	3	&	3	&	3.5	&	4	&	3\\
2572	&	5.5	&	5	&	4	&	5	&	4.5	&	4.5	&	5	&	6\\
4344	&	7	&	6	&	6	&	6.5	&	6	&	6.5	&	6.5	&	6\\
8883	&	9.5	&	8	&	7.5	&	7.5	&	7.5	&	9	&	8.5	&	8.5\\
21228	&	12	&	11.5	&	11	&	10	&	11.5	&	11	&	11.5	&	11\\
45808	&	14	&	13.5	&	13.5	&	14.5	&	12.5	&	13.5	&	13	&	21.5\\
92528	&	16.5	&	17.5	&	17	&	16.5	&	18	&	17.5	&	18.5	&	19.5\\
176179	&	24.5	&	24.5	&	22.5	&	20	&	23.5	&	20	&	23.5	&	25.5\\
\hline
\end{tabular}
\caption{Test Problem~\ref{pro19}:  Number of iterations of BiCGSTAB for uniform interior grids and 
boundary nodes from optimized triangulations as in Table~\ref{pr19gridfem}.
}
\label{pr19gridfem_iter}
\end{table}

In the last experiment  we generate interior nodes of $\Xi$ from the Halton sequence of quasi-random numbers
and the boundary nodes by the orthogonal projection. Surprisingly, the results are much better than  in the case of interior
grids combined with projected boundary nodes, see Table~\ref{pr19project} and Figure~\ref{Pr2_tables}(d). The parameters
of \texttt{oct-dist} are the same as in the case of Halton interior nodes for Test Problem~\ref{pro11}:
$\delta=0.9$, $s=1$, $n=3$ and $k=17$. We see that the best errors delivered by \texttt{oct-dist} and \texttt{pQR4sel} 
are not much higher than those achieved on the optimized triangulation. Selection methods \texttt{tet} and \texttt{oct}
failed to provide influence sets satisfying \eqref{exap}. The errors of \texttt{pQR3} and  \texttt{pQR4} 
are quite high, therefore we did not include them in the figure. The irregularities of the plots in Figure~\ref{Pr2_tables}(d) match well the behavior  of the
stability constant, see Table~\ref{pr19project_sigma}. Indeed, 
the pick in the plot of \texttt{oct-dist} when $i=1$ corresponds to $\sigma=546.9$, while otherwise $61.0\le\sigma\le80.1$ 
for this selection method. For \texttt{pQR4sel} we also see a pick (when $i=3$) corresponding to a relatively high
$\sigma=138.5$.  For \texttt{20near}, picks for $i=2$ and 6 also correlate with higher $\sigma$. Table~\ref{pr19project_sigma}
also include iteration numbers of BiCGSTAB that are in line with the observations made for previous types of nodes.

\begin{table}[htbp!]\scriptsize%
\centering
\renewcommand{\arraystretch}{1.25}
\begin{tabular}{|c||c|c|c||c|c|c|c||c|c|}
\hline
\multirow{2}{.7cm}{$\# {\Xi _{\rm int}}$} &\multicolumn{3}{c||}{Polyharmonic} & \multicolumn{2}{c|}{Pivoted QR}\\
\cline{2-6}
&	 {\scriptsize oct-dist} 
		& {\scriptsize 20near} & {\scriptsize pQR4sel}& {\scriptsize pQR3} & {\scriptsize pQR4}\\
\hline
 1165  & 6.9e{\tt-}02  & 9.2e{\tt-}02  & 1.0e{\tt-}01  & 9.2e{\tt-}01  & 3.3e{\tt+}00\\
 2474  & 2.7e{\tt-}01  & 7.7e{\tt-}02  & 6.4e{\tt-}02  & 8.0e{\tt-}02  & 3.7e{\tt-}01\\
 5187  & 2.7e{\tt-}02  & 8.3e{\tt-}02  & 5.2e{\tt-}02  & 2.9e{\tt+}00  & 8.1e{\tt-}01\\
 10754 & 1.9e{\tt-}02  & 5.3e{\tt-}02  & 5.0e{\tt-}02  & 2.3e{\tt-}02  & 2.5e{\tt-}01\\
 22093 & 1.3e{\tt-}02  & 1.5e{\tt-}02  & 1.5e{\tt-}02  & 3.3e{\tt-}02  & 3.6e{\tt-}01\\
 45172 & 7.6e{\tt-}03  & 1.8e{\tt-}02  & 8.1e{\tt-}03  & 2.0e{\tt-}02  & 2.8e{\tt-}01\\
 91893 & 6.0e{\tt-}03  & 2.8e{\tt-}02  & 6.5e{\tt-}03  & 1.6e{\tt-}02  & 3.7e{\tt-}02\\
 186274  & 4.0e{\tt-}03  & 9.0e{\tt-}03  & 4.2e{\tt-}03  & 5.1e{\tt-}02  & 5.5e{\tt-}02\\
\hline
density	&	16.3	&	18.8	&	18.9	&	9.6	&	18.9\\
\hline
\end{tabular}
\caption{Test Problem~\ref{pro19}: RRMS error $E_{\rm ref}$  for 
Halton interior nodes and projected boundary nodes.
}
\label{pr19project}
\end{table}

\begin{table}[htbp!]\scriptsize%
\centering
\renewcommand{\arraystretch}{1.25}
\begin{tabular}{|c||c|c||c|c||c|c||c|c||c|c|}
\hline
\multirow{2}{.7cm}{$\# {\Xi _{\rm int}}$} & \multicolumn{2}{|c||}{oct-dist}  & \multicolumn{2}{|c||}{20near} & \multicolumn{2}{|c||}{pQR4sel} 
& \multicolumn{2}{|c||}{pQR3} & \multicolumn{2}{|c|}{pQR4}\\
\cline{2-11} 
& $\sigma$ & $\#$iter & $\sigma$ & $\#$iter & $\sigma$ & $\#$iter & $\sigma$ & $\#$iter & $\sigma$ & $\#$iter\\
\hline
1165	&	64	&	4	&	1.9e{\tt +}03	&	5	&	582.7	&	4	&	3.4e{\tt +}03	&	9	&	4.3e{\tt +}05	&	232.5\\
2474	&	546.9	&	4.5	&	351	&	6.5	&	63	&	5.5	&	851.4	&	13	&	5.1e{\tt +}04	&	261.5\\
5187	&	62	&	5.5	&	245	&	8	&	61.8	&	10	&	7.5e{\tt +}03	&	9.5	&	5.3e{\tt +}04	&	316.5\\
10754	&	70	&	8	&	76.1	&	13	&	138.5	&	9	&	392.8	&	17.5	&	9.1e{\tt +}03	&	552.5\\
22093	&	61	&	10	&	63.0	&	16.5	&	61.5	&	12	&	3.8e{\tt +}16	&	fail	&	6.4e{\tt +}03	&	fail\\
45172	&	61	&	12	&	64.6	&	11.5	&	60.8	&	14	&	9.6e{\tt +}15	&	fail	&	5.2e{\tt +}03	&	fail\\
91893	&	80.2	&	15.5	&	118.4	&	23	&	61.3	&	24.5	&	1.2e{\tt +}16	&	fail	&	2.1e{\tt +}04	&	fail\\
186274	&	61.3	&	52.5	&	62.9	&	193.5	&	61.6	&	41	&	1.5e{\tt +}16	&	fail	&	4.2e{\tt +}04	&	fail \\
\hline
\end{tabular}
\caption{Test Problem~\ref{pro19}: stability constant $\sigma$ and number of BiCGSTAB  iterations $\#$iter for 
Halton interior nodes and projected boundary nodes as in Table~\ref{pr19project}. 
}
\label{pr19project_sigma}
\end{table}

\begin{figure}[!ht]
\begin{center}
\includegraphics[width=7.5cm,height=5.0cm]{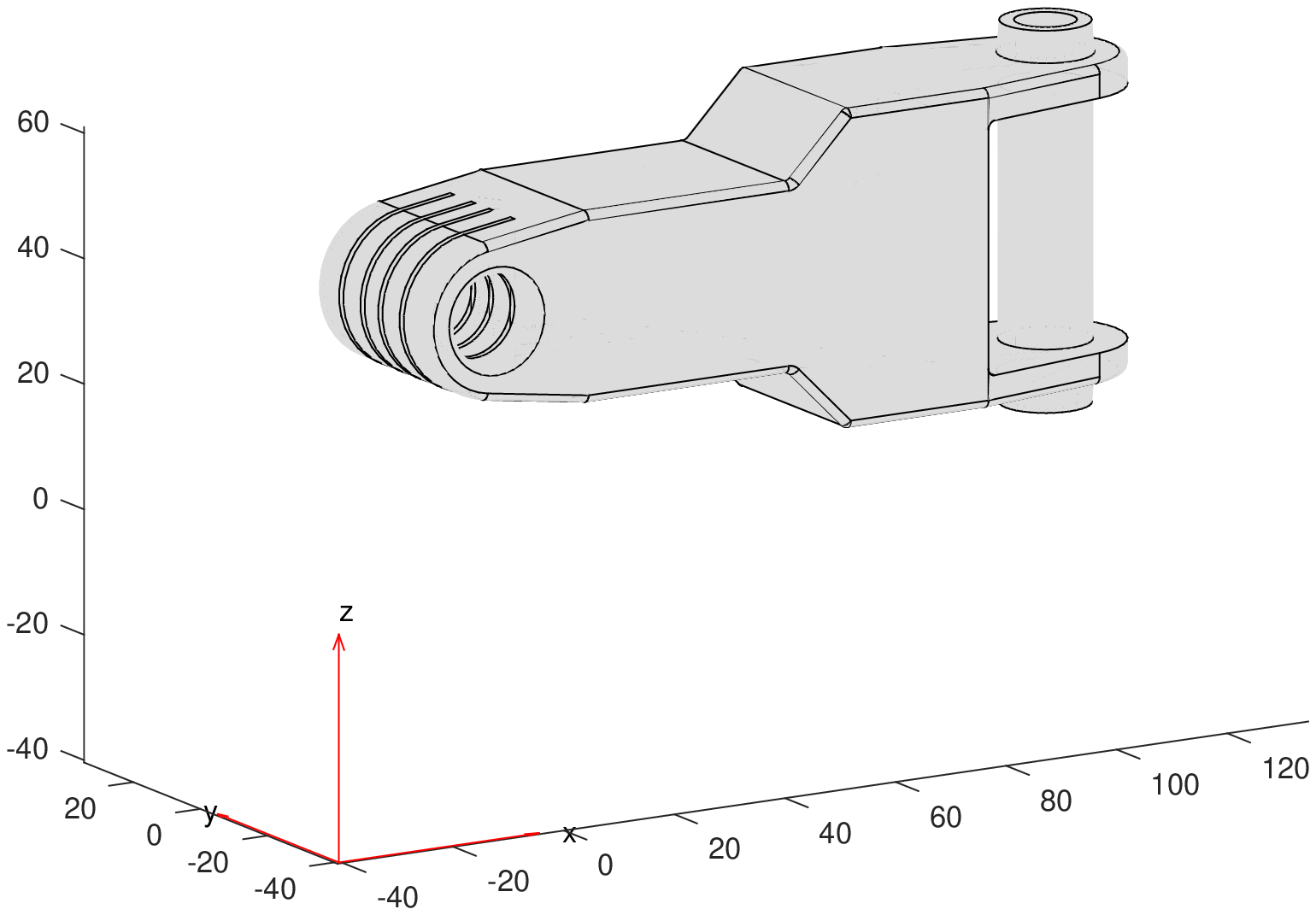}
\caption{Test Problem~\ref{pro5}: Domain \texttt{ForearmLink}.
\label{fig5_1}
}
\end{center}
\end{figure}

\begin{testproblem}[\texttt{ForearmLink}]
\label{pro5}
Poisson equation $\Delta u = -10$ with zero Dirichlet boundary conditions on the domain $\Omega$ defined in
the STL file  `ForearmLink.stl' shipped with MATLAB  PDE Toolbox \cite{PDEtool},
see Figure~\ref{fig5_1}.
\end{testproblem}

A reference solution is computed by  the second order finite element method on a triangulation with 641825 nodes 
obtained by \texttt{generateMesh} with
$\texttt{Hmax}=6.8/2^{7/3}\approx 1.35$, $\texttt{Hmin}=\texttt{Hmax}/3$,
\texttt{Hgrad}=1.5 and \texttt{GeometricOrder}=2. 

\begin{figure}[htbp!]%
\begin{center}
\hbox{\subfigure[Optimized triangulation]{\includegraphics[width=7.5cm,height=5.0cm]{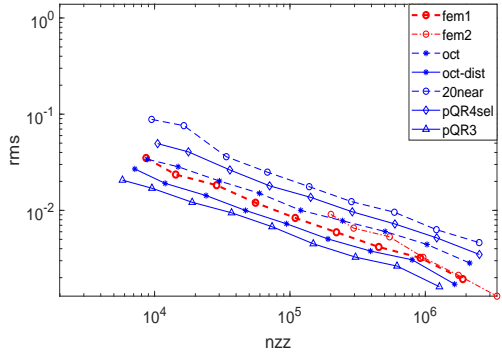}}
\subfigure[Unoptimized triangulation]{\includegraphics[width=7.5cm,height=5.0cm]{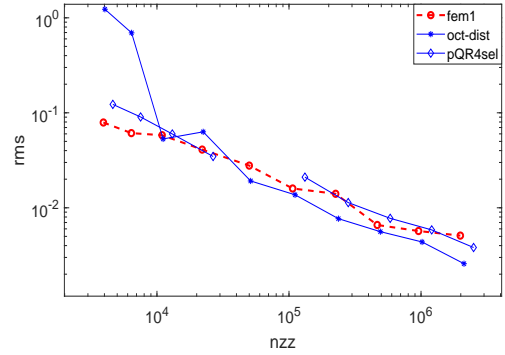}}}
\hbox{\subfigure[Interior grid and boundary nodes of (a)]{\includegraphics[width=7.5cm,height=5.0cm]{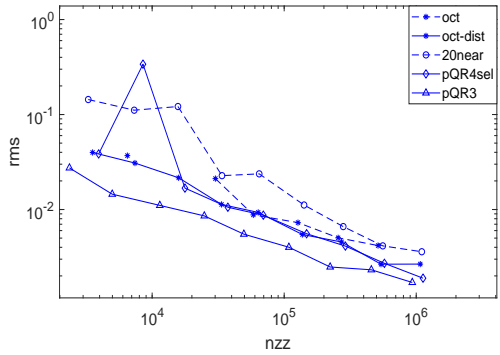}}
\subfigure[Halton interior and projected boundary nodes]{\includegraphics[width=7.5cm,height=5.0cm]{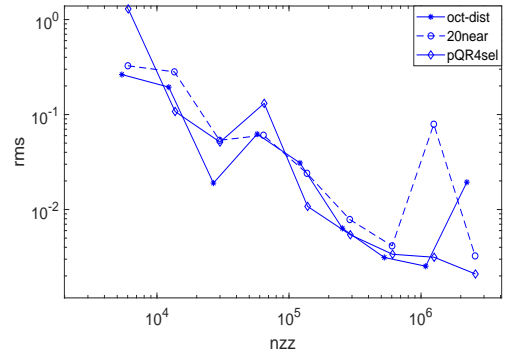}}}
\caption{Test Problem~\ref{pro5}: RMS errors as functions of \texttt{nnz} for meshless FD method 
with various versions of %
stencil selection, and finite element method. 
(a)  Tables~\ref{pr5fem} and \ref{pr5fem2}; (b) Table~\ref{pr5gmsh}; (c) Table~\ref{pr5gridfem}; (d) Table~\ref{pr5project}.
\label{Pr3_tables}
}
\end{center}
\end{figure}

Similar to the previous test problems, we run several numerical experiments using nodes of different types:
vertices of optimized/unoptimized triangulations, interior Cartesian grids combined with boundary vertices of the
optimized triangulation, and  Halton interior nodes with projected nodes on the boundary. 
The results are presented in Tables~\ref{pr5fem}--\ref{pr5project} and Figure~\ref{Pr3_tables}. For \texttt{oct-dist}
we use $\delta=0.7$, with other parameters $s,n,k$ reported in the captions of respective tables. 
Note that we omitted  the high error of \texttt{pQR4sel} for $i=4$ in its plot in Figure~\ref{Pr3_tables}(b) as well as
several parts of the plot for \texttt{oct} in Figure~\ref{Pr3_tables}(c) with \texttt{NaN} in the respective table.

The results are similar to those for Test Problem~\ref{pro19}: \texttt{tet} fails to satisfy \eqref{exap}, 
and the same often happens to \texttt{20near} and \texttt{oct} on less regular node sets;  
the error plots in Figure~\ref{Pr3_tables}(b-d) are not always monotone, which indicates fluctuations in the stability
of the system matrix; \texttt{oct-dist} belongs to best performers overall in all cases, and \texttt{pQR3} 
for optimized triangulations and gridded nodes. 
In contrast to Test Problem~\ref{pro19}, 
\texttt{pQR4} does not perform well on the unoptimized triangulations until they are sufficiently fine. Note that the second
order finite element method \texttt{fem2} does not seem advantageous over \texttt{fem1} for this domain, especially 
if the higher density of its
system matrix is taken into account. Nevertheless, the slope of its error plot in Figure~\ref{Pr3_tables}(a)  indicates a
higher convergence order, which presumably would make \texttt{fem2} preferable on further refined triangulations.

\begin{table}[htbp!]\scriptsize %
\centering
\renewcommand{\arraystretch}{1.25}
\begin{tabular}{|c||c||c|c|c|c|c||c|c|}
\hline
\multirow{2}{.7cm}{$\# {\Xi _{\rm int}}$} & FEM   
 &	\multicolumn{5}{|c||}{Polyharmonic $r^5$ with quadratic polynomial} & \multicolumn{2}{c|}{Pivoted QR}\\
\cline{2-9}
		& {\scriptsize fem1}	& {\scriptsize tet} &	{\scriptsize oct} & {\scriptsize oct-dist} 
		& {\scriptsize 20near} & {\scriptsize pQR4sel} & {\scriptsize pQR3} & {\scriptsize pQR4}\\
\hline
816	&	3.5e{\tt-}02	&	NaN	&	3.4e{\tt-}02	&	2.7e{\tt-}02	&	8.8e{\tt-}02	&	4.9e{\tt-}02	&	2.1e{\tt-}02	&	3.0e{\tt+}00\\
1254	&	2.4e{\tt-}02	&	NaN	&	2.8e{\tt-}02	&	1.9e{\tt-}02	&	7.6e{\tt-}02	&	4.1e{\tt-}02	&	1.7e{\tt-}02	&	9.4e{\tt-}02\\
2414	&	1.8e{\tt-}02	&	NaN	&	2.0e{\tt-}02	&	1.4e{\tt-}02	&	3.6e{\tt-}02	&	2.6e{\tt-}02	&	1.2e{\tt-}02	&	4.4e{\tt-}02\\
4423	&	1.2e{\tt-}02	&	NaN	&	1.5e{\tt-}02	&	1.0e{\tt-}02	&	2.5e{\tt-}02	&	1.8e{\tt-}02	&	9.4e{\tt-}03	&	4.2e{\tt-}02\\
8401	&	8.4e{\tt-}03	&	NaN	&	1.0e{\tt-}02	&	7.3e{\tt-}03	&	1.8e{\tt-}02	&	1.4e{\tt-}02	&	6.8e{\tt-}03	&	6.9e{\tt-}03\\
16437	&	5.9e{\tt-}03	&	NaN	&	7.8e{\tt-}03	&	5.1e{\tt-}03	&	1.2e{\tt-}02	&	9.7e{\tt-}03	&	4.5e{\tt-}03	&	1.0e{\tt-}02\\
33030	&	4.2e{\tt-}03	&	NaN	&	6.0e{\tt-}03	&	3.8e{\tt-}03	&	9.6e{\tt-}03	&	7.2e{\tt-}03	&	3.3e{\tt-}03	&	2.7e{\tt-}03\\
65652	&	3.2e{\tt-}03	&	NaN	&	4.4e{\tt-}03	&	3.1e{\tt-}03	&	6.3e{\tt-}03	&	5.2e{\tt-}03	&	2.6e{\tt-}03	&	3.5e{\tt-}03\\
133295	&	1.9e{\tt-}03	&	NaN	&	2.8e{\tt-}03	&	1.7e{\tt-}03	&	4.6e{\tt-}03	&	3.5e{\tt-}03	&	1.6e{\tt-}03	&	1.3e{\tt-}03\\
\hline
density	&	14.2	&	14.2	&	16.0	&	12.4	&	18.8	&	18.8	&	9.5	&	18.8\\
\hline
\end{tabular}
\caption{Test Problem~\ref{pro5}:  RRMS errors $E_{\rm ref}$ for optimized
triangulations with $H_0=6.8$. Parameters of \texttt{oct-dist}: $s=1$, $n=3$,  $k=13$.
}
\label{pr5fem}
\end{table}

\begin{table}[htbp!]\scriptsize %
\centering
\renewcommand{\arraystretch}{1.25}
\begin{tabular}{|c|c|c|c|c|c|c|c|c||c|}
\hline
$\# {\Xi _{\rm int}}$ & 10268 & 	14078 & 	24426 & 	40806 & 	71249 & 	132581 & density\\
\hline
\hline
fem2	&	9.1e{\tt-}03	&	6.5e{\tt-}03	&	5.3e{\tt-}03	&	3.3e{\tt-}03	&	2.1e{\tt-}03	&	1.3e{\tt-}03	&	25.5 \\
\hline
\end{tabular}
\caption{Test Problem~\ref{pro5}:   RRMS errors $E_{\rm ref}$ of the quadratic finite element method 
for optimized triangulations with $H_0=7.1$.}
\label{pr5fem2}
\end{table}

\begin{table}[htbp!]\scriptsize %
\centering
\renewcommand{\arraystretch}{1.25}
\begin{tabular}{|c||c||c|c|c|c|c||c|c|}
\hline
\multirow{2}{.7cm}{$\# {\Xi _{\rm int}}$} & FEM   
 &	\multicolumn{5}{|c||}{Polyharmonic $r^5$ with quadratic polynomial} & \multicolumn{2}{c|}{Pivoted QR}\\
\cline{2-9}
		& {\scriptsize fem1}	& {\scriptsize tet} &	{\scriptsize oct} & {\scriptsize oct-dist} 
		& {\scriptsize 20near} & {\scriptsize pQR4sel} & {\scriptsize pQR3} & {\scriptsize pQR4}\\
\hline
553	&	7.9e{\tt-}02	&	NaN	&	NaN	&	1.23e{\tt+}00	&	NaN	&	1.23e{\tt-}01	&	1.85e{\tt+}00	&	3.04e{\tt+}00\\
738	&	6.1e{\tt-}02	&	NaN	&	NaN	&	6.95e{\tt-}01	&	NaN	&	9.04e{\tt-}02	&	2.22e{\tt-}01	&	1.79e{\tt+}00\\
1092	&	5.8e{\tt-}02	&	NaN	&	NaN	&	5.31e{\tt-}02	&	NaN	&	5.98e{\tt-}02	&	1.30e{\tt-}01	&	3.62e{\tt-}01\\
1918	&	4.1e{\tt-}02	&	NaN	&	NaN	&	6.30e{\tt-}02	&	NaN	&	3.47e{\tt-}02	&	1.14e{\tt-}01	&	1.04e{\tt+}01\\
3967	&	2.8e{\tt-}02	&	NaN	&	NaN	&	1.92e{\tt-}02	&	NaN	&	1.90e{\tt+}05	&	1.06e{\tt+}01	&	6.02e{\tt-}01\\
7964	&	1.6e{\tt-}02	&	NaN	&	NaN	&	1.37e{\tt-}02	&	NaN	&	2.10e{\tt-}02	&	1.20e{\tt-}01	&	5.89e{\tt-}02\\
16219	&	1.4e{\tt-}02	&	NaN	&	NaN	&	7.68e{\tt-}03	&	NaN	&	1.13e{\tt-}02	&	Inf	&	8.26e{\tt-}03\\
32460	&	6.6e{\tt-}03	&	NaN	&	NaN	&	5.61e{\tt-}03	&	NaN	&	7.74e{\tt-}03	&	2.27e{\tt-}01	&	1.00e{\tt-}01\\
65240	&	5.7e{\tt-}03	&	NaN	&	NaN	&	4.37e{\tt-}03	&	NaN	&	5.85e{\tt-}03	&	2.46e{\tt+}06	&	2.25e{\tt-}02\\
132782	&	5.1e{\tt-}03	&	NaN	&	NaN	&	2.58e{\tt-}03	&	NaN	&	3.82e{\tt-}03	&	3.63e{\tt-}03	&	1.76e{\tt-}03\\
\hline
density	&	14.9	&	14.9	&	15.9	&	15.93	&	18.3	&	18.8	&	9.6	&	18.8\\
\hline
\end{tabular}
\caption{Test Problem~\ref{pro5}:  RRMS errors $E_{\rm ref}$ for unoptimized triangulations.
Parameters of \texttt{oct-dist}: $s=3$, $n=6$ and  $k=17$.
}
\label{pr5gmsh}
\end{table}

\begin{table}[htbp!]\scriptsize%
\centering
\renewcommand{\arraystretch}{1.25}
\begin{tabular}{|c||c|c|c|c|c|c|}
\hline
Triangulation	&	$\min\gamma$&	$\avg\gamma$	&	$0<\gamma\leqslant0.25$	&	$0.25<\gamma\leqslant0.5$	
              &	$0.5<\gamma\leqslant0.75$	&	$0.75<\gamma\leqslant1.0$	\\
\hline
Optimized	&	0.09& 0.88	&	0.0\%  	&	0.1\% &	4.6\%	&	95.3\%\\
Unoptimized	&	1.0e{\tt-}13& 0.72	&	7.7\%	&	8.7\% &	24.6\%	&	58.9\%	 \\
\hline
\end{tabular}
\caption{Test Problem~\ref{pro5}: Statistics of the aspect ratio $\gamma_T$ for the simplices  of the  3D 
triangulations used in the tests with finite element method \texttt{fem1} reported in Tables~\ref{pr5fem} and 
\ref{pr5gmsh}.
}
\label{pr5meshq}
\end{table}

\begin{table}[htbp!]\scriptsize %
\centering
\renewcommand{\arraystretch}{1.25}
\begin{tabular}{|c|c|c|c|c||c|c|}
\hline
\multirow{2}{.7cm}{$\# {\Xi _{\rm int}}$} &\multicolumn{4}{c||}{Polyharmonic $r^5$ with quadratic polynomial} & \multicolumn{2}{c|}{Pivoted QR}\\
\cline{2-7}
& {\scriptsize oct} & {\scriptsize oct-dist} & {\scriptsize 20near} & {\scriptsize pQR4sel} & {\scriptsize pQR3} & {\scriptsize pQR4}\\
\hline
393	&	NaN	&	4.0e{\tt-}02	&	1.4e{\tt-}01	&	3.9e{\tt-}02	&	2.7e{\tt-}02	&	3.1e{\tt-}01\\
805	&	3.7e{\tt-}02	&	3.1e{\tt-}02	&	1.1e{\tt-}01	&	3.4e{\tt-}01	&	1.5e{\tt-}02	&	8.4e{\tt-}02\\
1812	&	NaN	&	2.2e{\tt-}02	&	1.2e{\tt-}01	&	1.7e{\tt-}02	&	1.1e{\tt-}02	&	2.7e{\tt-}02\\
3816	&	2.1e{\tt-}02	&	1.1e{\tt-}02	&	2.3e{\tt-}02	&	1.1e{\tt-}02	&	8.6e{\tt-}03	&	6.7e{\tt-}02\\
7507	&	8.8e{\tt-}03	&	9.3e{\tt-}03	&	2.4e{\tt-}02	&	8.7e{\tt-}03	&	5.5e{\tt-}03	&	2.1e{\tt-}02\\
16143	&	7.3e{\tt-}03	&	5.5e{\tt-}03	&	1.1e{\tt-}02	&	5.6e{\tt-}03	&	4.0e{\tt-}03	&	1.8e{\tt-}02\\
32862	&	5.0e{\tt-}03	&	4.6e{\tt-}03	&	6.6e{\tt-}03	&	4.1e{\tt-}03	&	2.5e{\tt-}03	&	6.2e{\tt-}03\\
66716	&	4.2e{\tt-}03	&	2.7e{\tt-}03	&	4.1e{\tt-}03	&	2.7e{\tt-}03	&	2.3e{\tt-}03	&	9.1e{\tt-}03\\
135410	&	NaN	&	2.7e{\tt-}03	&	3.6e{\tt-}03	&	1.9e{\tt-}03	&	1.7e{\tt-}03	&	7.7e{\tt+}00\\
\hline
density	&	7.6	&	7.9	&	8.1	&	8.3	&	6.9	&	8.3\\
\hline
\end{tabular}
\caption{Test Problem~\ref{pro5}: RRMS error  $E_{\rm ref}$ for uniform interior grids and 
boundary nodes from optimized triangulations. The 7-node grid stencil \eqref{7star} is used whenever possible.
Parameters of \texttt{oct-dist}: $s=1$, $n=3$ and  $k=18$.
}
\label{pr5gridfem}
\end{table}

\begin{table}[htbp!]\scriptsize%
\centering
\renewcommand{\arraystretch}{1.25}
\begin{tabular}{|c||c|c|c||c|c|c|c||c|c|}
\hline
\multirow{2}{.7cm}{$\# {\Xi _{\rm int}}$} &\multicolumn{3}{c||}{Polyharmonic} & \multicolumn{2}{c|}{Pivoted QR}\\
\cline{2-6}
&	 {\scriptsize oct-dist}  & {\scriptsize 20near} & {\scriptsize pQR4sel}& {\scriptsize pQR3} & {\scriptsize pQR4}\\
\hline
 392	&	2.6e-01	&	3.3e-01	&	1.3e+00	&	1.0e-01	&	3.3e-01\\
862	&	1.9e-01	&	2.8e-01	&	1.1e-01	&	5.3e-01	&	3.7e-01\\
1819	&	1.9e-02	&	5.4e-02	&	5.1e-02	&	2.0e-01	&	1.5e+00\\
3808	&	6.2e-02	&	6.1e-02	&	1.3e-01	&	6.2e-02	&	8.9e-01\\
7886	&	3.1e-02	&	2.4e-02	&	1.1e-02	&	2.6e-02	&	1.6e+00\\
16181	&	6.3e-03	&	7.9e-03	&	5.4e-03	&	2.2e-02	&	2.9e-01\\
33076	&	3.1e-03	&	4.1e-03	&	3.4e-03	&	3.9e-03	&	3.4e-02\\
67148	&	2.5e-03	&	7.9e-02	&	3.1e-03	&	3.7e-02	&	3.5e-01\\
135971	&	1.9e-02	&	3.3e-03	&	2.1e-03	&	2.1e+00	&	2.1e-02\\
\hline
density	&	16.4	&	18.9	&	19.0	&	9.6	&	19.0\\
\hline
\end{tabular}
\caption{
Test Problem~\ref{pro5}: RRMS error $E_{\rm ref}$  for 
Halton interior nodes and projected boundary nodes.
Parameters of \texttt{oct-dist}: $s=1$, $n=3$ and  $k=17$.
}
\label{pr5project}
\end{table}

\begin{figure}[htbp!]%
\begin{center}
\includegraphics[width=7.5cm,height=5.0cm]{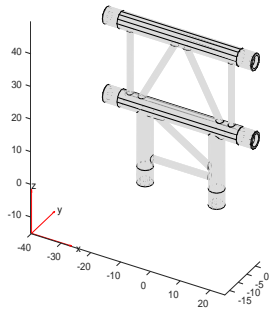}
\caption{Test Problem~\ref{pro13}: Domain \texttt{BeamTrussJunction}  visualized 
by MATLAB PDE Toolbox command {\tt pdegplot}.
\label{fig13_1}
}
\end{center}
\end{figure}

\begin{testproblem}[\texttt{BeamTrussJunction}] 
\label{pro13}
Poisson equation $\Delta u = -10$ with zero Dirichlet boundary conditions on the domain $\Omega$ 
defined in the STL model ``Beam Truss Cross and T Junction 134'' designed by akerStudio,
 available from the platform Cults at the link 
{\rm\url{ https://cults3d.com/en/3d-model/various/beam-truss-cross-and-t-junction-134}} 
The domain is visualized in Figure~\ref{fig13_1}.
\end{testproblem}

A reference solution is computed by  the second order finite element method on a triangulation with 1337199 nodes 
obtained by \texttt{generateMesh} with
$\texttt{Hmax}= 0.3$, $\texttt{Hmin}=\texttt{Hmax}/3$,
\texttt{Hgrad}=1.5 and \texttt{GeometricOrder}=2. 

\begin{figure}[!ht]
\begin{center}
\hbox{\subfigure[Optimized triangulation]{\includegraphics[width=7.5cm,height=5.0cm]{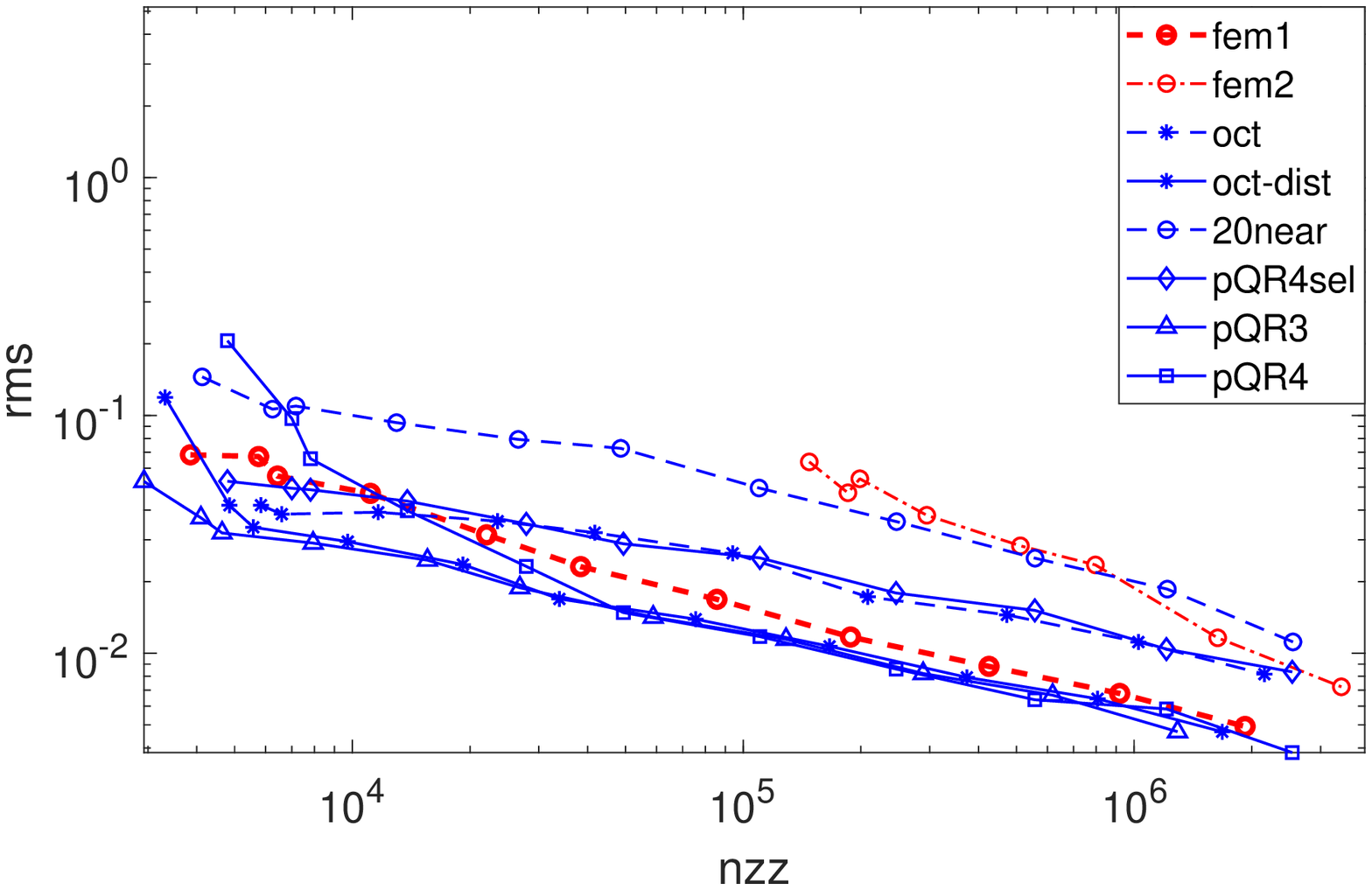}}
\subfigure[Unoptimized triangulation]{\includegraphics[width=7.5cm,height=5.0cm]{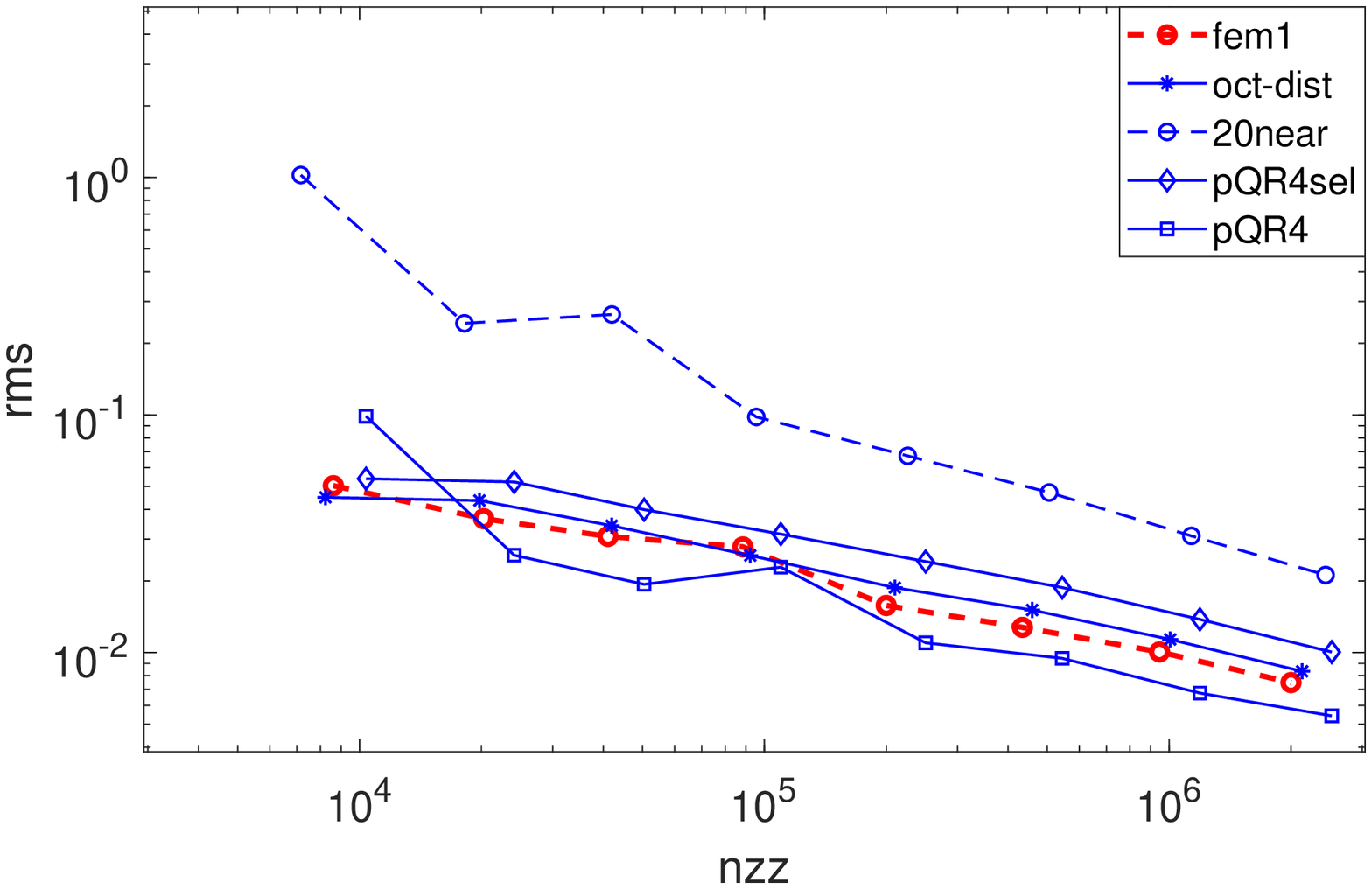}}}
\hbox{\subfigure[Interior grid and boundary nodes of (a)]{\includegraphics[width=7.5cm,height=5.0cm]{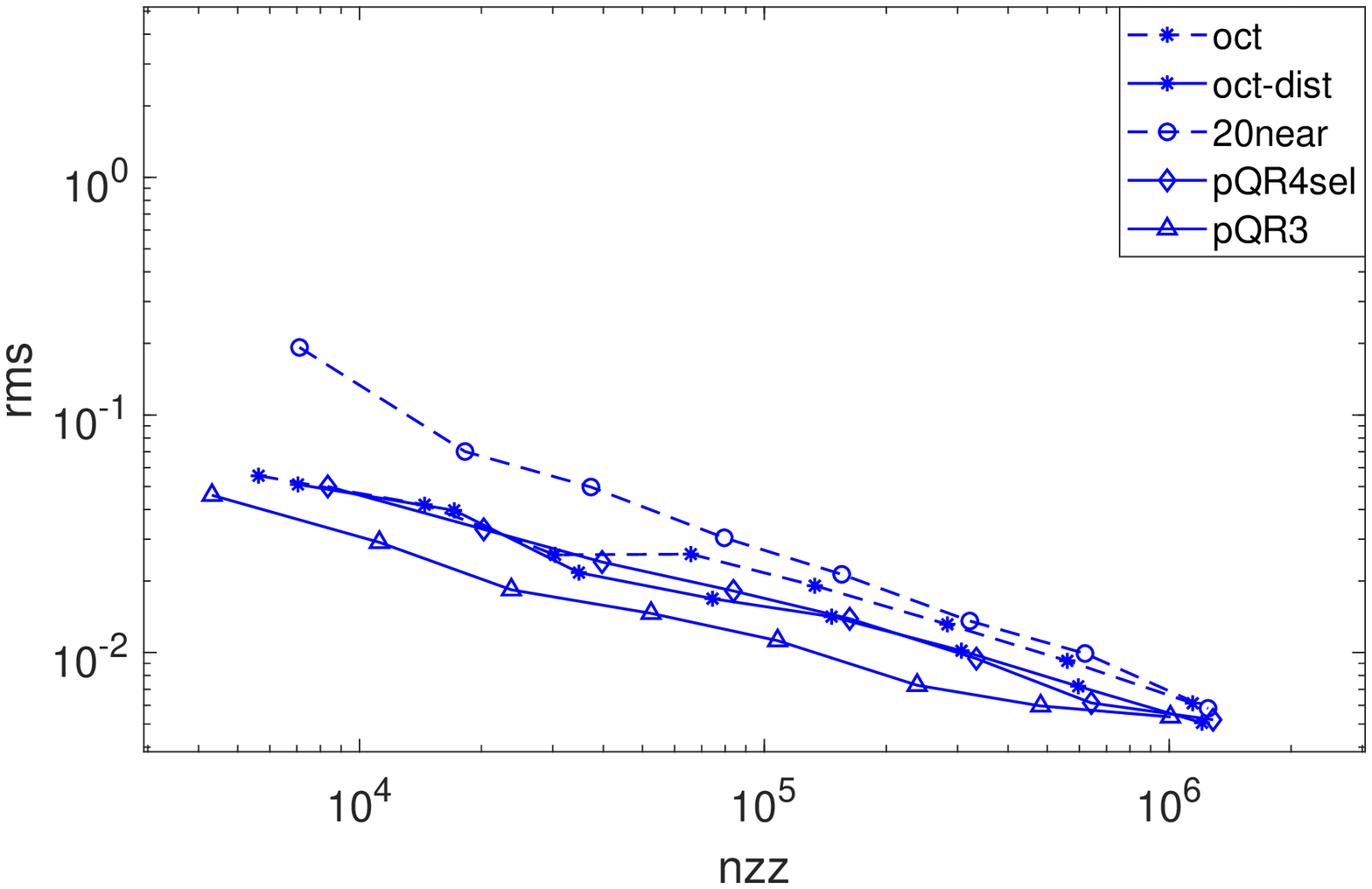}}
\subfigure[Halton interior and projected boundary nodes]{\includegraphics[width=7.5cm,height=5.0cm]{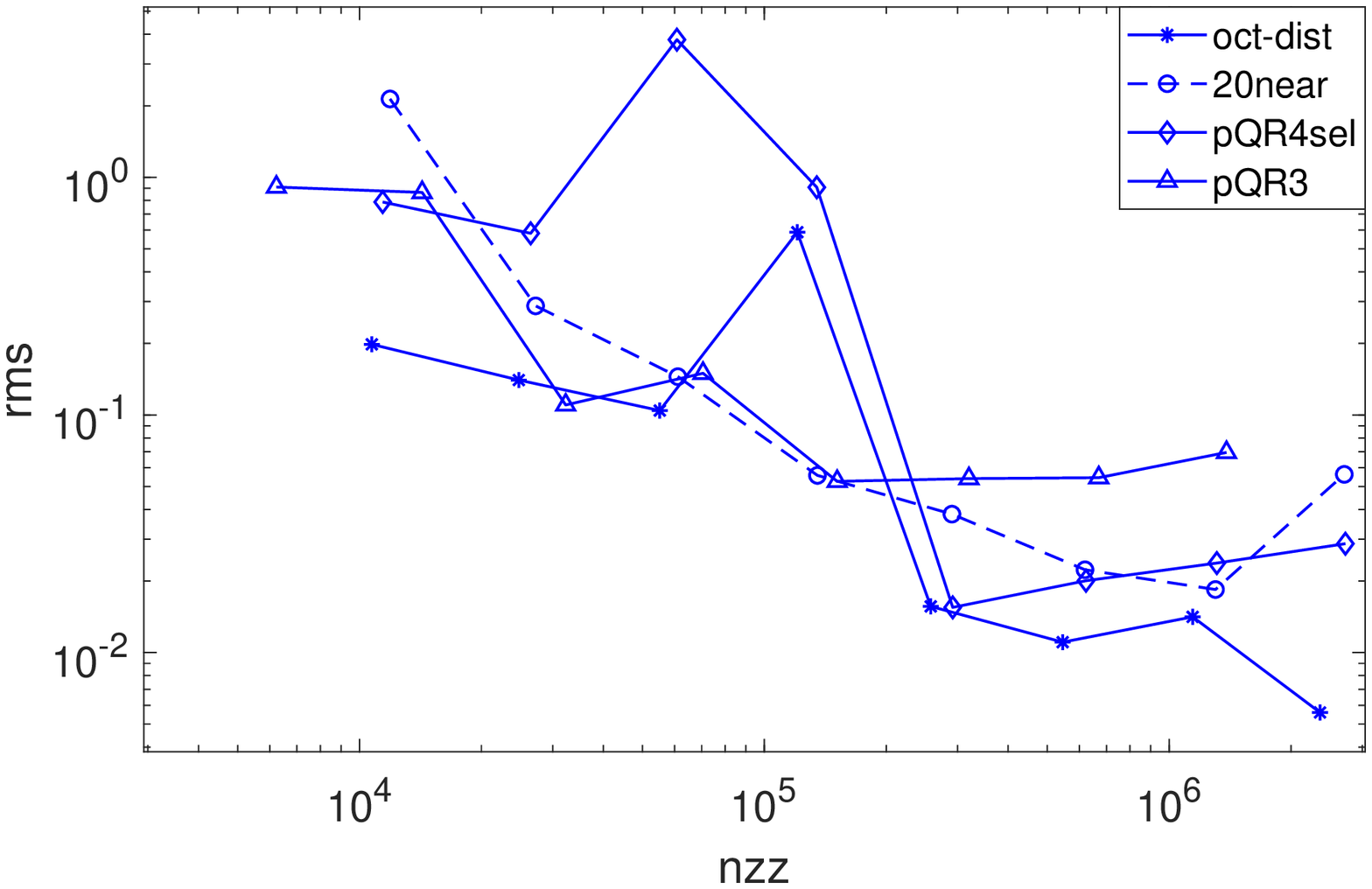}}}
\caption{Test Problem~\ref{pro13}: RMS errors as functions of \texttt{nnz} for meshless FD method 
with various versions of %
stencil selection, and finite element method. 
(a)  Tables~\ref{pr13fem} and \ref{pr13fem2}; (b) Table~\ref{pr13gmsh}; (c) Table~\ref{pr13gridfem}; (d) Table~\ref{pr13project}
\label{fig13tables}
}
\end{center}
\end{figure}

We run the same type of experiments as for Test Problem~\ref{pro5}. The results are reported in 
Tables~\ref{pr13fem}--\ref{pr13project} and Figure~\ref{fig13tables}. We use $\delta=0.7$ for \texttt{oct-dist}, 
with other parameters $s,n,k$ provided in the captions of respective tables. The results confirm the same observations made on
Test Problems~\ref{pro19} and \ref{pro5}. Remarkably, the errors of \texttt{fem2} are now significantly higher than those
of  \texttt{fem1} and of several mFD methods despite its higher density of the system matrix. We still expect it would 
outperform them on denser node sets thanks to its higher order. However, assuming that the accuracy seen in the tables, e.g.\ the
RRMS error of $10^{-2}$, is sufficient for an application, there seems to be no point in using a higher order method that
achieves a desired error threshold with many more degrees of freedom.

\begin{table}[htbp!]\scriptsize %
\centering
\renewcommand{\arraystretch}{1.25}
\begin{tabular}{|c||c||c|c|c|c|c||c|c|}
\hline
\multirow{2}{.7cm}{$\# {\Xi _{\rm int}}$} & FEM   
 &	\multicolumn{5}{|c||}{Polyharmonic $r^5$ with quadratic polynomial} & \multicolumn{2}{c|}{Pivoted QR}\\
\cline{2-9}
		& {\scriptsize fem1}	& {\scriptsize tet} &	{\scriptsize oct} & {\scriptsize oct-dist} 
		& {\scriptsize 20near} & {\scriptsize pQR4sel} & {\scriptsize pQR3} & {\scriptsize pQR4}\\
\hline
538	&	6.8e{\tt-}02	&	NaN	&	NaN	&	1.2e{\tt-}01	&	1.5e{\tt-}01	&	5.3e{\tt-}02	&	5.3e{\tt-}02	&	2.1e{\tt-}01\\
664	&	6.7e{\tt-}02	&	NaN	&	4.2e{\tt-}02	&	4.2e{\tt-}02	&	1.1e{\tt-}01	&	4.9e{\tt-}02	&	3.7e{\tt-}02	&	9.7e{\tt-}02\\
760	&	5.6e{\tt-}02	&	NaN	&	3.8e{\tt-}02	&	3.4e{\tt-}02	&	1.1e{\tt-}01	&	4.9e{\tt-}02	&	3.2e{\tt-}02	&	6.6e{\tt-}02\\
1177	&	4.7e{\tt-}02	&	2.9e{\tt-}02	&	3.9e{\tt-}02	&	3.0e{\tt-}02	&	9.3e{\tt-}02	&	4.4e{\tt-}02	&	2.9e{\tt-}02	&	4.0e{\tt-}02\\
2122	&	3.2e{\tt-}02	&	2.1e{\tt-}02	&	3.6e{\tt-}02	&	2.4e{\tt-}02	&	7.9e{\tt-}02	&	3.5e{\tt-}02	&	2.5e{\tt-}02	&	2.3e{\tt-}02\\
3421	&	2.3e{\tt-}02	&	1.7e{\tt-}02	&	3.2e{\tt-}02	&	1.7e{\tt-}02	&	7.3e{\tt-}02	&	2.9e{\tt-}02	&	1.9e{\tt-}02	&	1.5e{\tt-}02\\
7101	&	1.7e{\tt-}02	&	NaN	&	2.6e{\tt-}02	&	1.4e{\tt-}02	&	5.0e{\tt-}02	&	2.5e{\tt-}02	&	1.4e{\tt-}02	&	1.2e{\tt-}02\\
14798	&	1.2e{\tt-}02	&	1.1e{\tt-}02	&	1.7e{\tt-}02	&	1.1e{\tt-}02	&	3.6e{\tt-}02	&	1.8e{\tt-}02	&	1.1e{\tt-}02	&	8.6e{\tt-}03\\
32246	&	8.8e{\tt-}03	&	NaN	&	1.5e{\tt-}02	&	7.9e{\tt-}03	&	2.5e{\tt-}02	&	1.5e{\tt-}02	&	8.2e{\tt-}03	&	6.4e{\tt-}03\\
67566	&	6.8e{\tt-}03	&	NaN	&	1.1e{\tt-}02	&	6.4e{\tt-}03	&	1.9e{\tt-}02	&	1.0e{\tt-}02	&	6.7e{\tt-}03	&	5.8e{\tt-}03\\
138312	&	4.9e{\tt-}03	&	NaN	&	8.2e{\tt-}03	&	4.7e{\tt-}03	&	1.1e{\tt-}02	&	8.4e{\tt-}03	&	4.7e{\tt-}03	&	3.8e{\tt-}03\\
\hline
density	&	13.9	&	13.9	&	15.6	&	12.2	&	18.4	&	18.4	&	9.3	&	18.4\\
\hline
\end{tabular}
\caption{Test Problem~\ref{pro13}:  RRMS errors $E_{\rm ref}$ for optimized
triangulations with $H_0=3.0$. Parameters of \texttt{oct-dist}: $s=1$, $n=3$,  $k=13$.
}
\label{pr13fem}
\end{table}

\begin{table}[htbp!]\scriptsize %
\centering
\renewcommand{\arraystretch}{1.25}
\begin{tabular}{|c|c|c|c|c|c|c|c|c||c|}
\hline
$\# {\Xi _{\rm int}}$ &8809&	10518&	11200&	15396&	24735&	35879&	69018&	135949 & density\\
\hline
\hline
fem2						& 6.4e{\tt-}02	&	4.7e{\tt-}02	&	5.4e{\tt-}02	&	3.8e{\tt-}02	&	2.8e{\tt-}02	&	2.4e{\tt-}02	&	1.2e{\tt-}02	&7.2e{\tt-}03&	24.9  \\
\hline
\end{tabular}
\caption{Test Problem~\ref{pro13}:   RRMS errors $E_{\rm ref}$ of the quadratic finite element method 
for optimized triangulations with $H_0=3.0$. 
}
\label{pr13fem2}
\end{table}

\begin{table}[htbp!]\scriptsize %
\centering
\renewcommand{\arraystretch}{1.25}
\begin{tabular}{|c||c||c|c|c|c|c||c|c|}
\hline
\multirow{2}{.7cm}{$\# {\Xi _{\rm int}}$} & FEM   
 &	\multicolumn{5}{|c||}{Polyharmonic $r^5$ with quadratic polynomial} & \multicolumn{2}{c|}{Pivoted QR}\\
\cline{2-9}
		& {\scriptsize fem1}	& {\scriptsize tet} &	{\scriptsize oct} & {\scriptsize oct-dist} 
		& {\scriptsize 20near} & {\scriptsize pQR4sel} & {\scriptsize pQR3} & {\scriptsize pQR4}\\
\hline
876	&	5.0e{\tt-}02	&	NaN	&	NaN	&	4.5e{\tt-}02	&	1.0e{\tt+}00	&	5.4e{\tt-}02	&	1.5e{\tt-}01	&	9.9e{\tt-}02\\
1852	&	3.7e{\tt-}02	&	NaN	&	NaN	&	4.4e{\tt-}02	&	2.4e{\tt-}01	&	5.2e{\tt-}02	&	2.6e{\tt+}01	&	2.6e{\tt-}02\\
3483	&	3.1e{\tt-}02	&	NaN	&	NaN	&	3.4e{\tt-}02	&	2.6e{\tt-}01	&	4.0e{\tt-}02	&	4.4e{\tt-}01	&	1.9e{\tt-}02\\
7137	&	2.8e{\tt-}02	&	NaN	&	3.3e{\tt-}02	&	2.6e{\tt-}02	&	9.8e{\tt-}02	&	3.2e{\tt-}02	&	3.5e{\tt-}02	&	2.3e{\tt-}02\\
15226	&	1.6e{\tt-}02	&	NaN	&	2.3e{\tt-}02	&	1.9e{\tt-}02	&	6.7e{\tt-}02	&	2.4e{\tt-}02	&	1.5e{\tt-}02	&	1.1e{\tt-}02\\
31609	&	1.3e{\tt-}02	&	2.9e{\tt-}02	&	1.7e{\tt-}02	&	1.5e{\tt-}02	&	4.7e{\tt-}02	&	1.9e{\tt-}02	&	9.7e{\tt-}03	&	9.5e{\tt-}03\\
66627	&	1.0e{\tt-}02	&	NaN	&	NaN	&	1.1e{\tt-}02	&	3.1e{\tt-}02	&	1.4e{\tt-}02	&	7.4e{\tt-}03	&	6.8e{\tt-}03\\
137521	&	7.5e{\tt-}03	&	NaN	&	9.2e{\tt-}03	&	8.3e{\tt-}03	&	2.1e{\tt-}02	&	1.0e{\tt-}02	&	6.0e{\tt-}03	&	5.4e{\tt-}03\\
\hline
density	&	14.6	&	14.6	&	15.5	&	15.5	&	17.7	&	18.3	&	9.4	&	18.3\\
\hline
\end{tabular}
\caption{ Test Problem~\ref{pro13}:  RRMS errors $E_{\rm ref}$ for unoptimized triangulations 
with $H_0=1.5/2^{1/3}$. %
Parameters of \texttt{oct-dist}: $s=3$, $n=6$ and  $k=17$.
}
\label{pr13gmsh}
\end{table}

\begin{table}[htbp!]\scriptsize%
\centering
\renewcommand{\arraystretch}{1.25}
\begin{tabular}{|c||c|c|c|c|c|c|}
\hline
Triangulation	&	$\min\gamma$&	$\avg\gamma$	&	$0<\gamma\leqslant0.25$	&	$0.25<\gamma\leqslant0.5$	
              &	$0.5<\gamma\leqslant0.75$	&	$0.75<\gamma\leqslant1.0$	\\
\hline
Optimized	&	0.43& 0.88	&	0.0\%  	&	0.0\% &	6.0\%	&	94.0\%\\
Unoptimized	&	1.2e{\tt-}05& 0.75	&	3.6\%	&	7.3\% &	28.1\%	&	61.0\%	 \\
\hline
\end{tabular}
\caption{Test Problem~\ref{pro13}:  Statistics of the aspect ratio $\gamma_T$ for the simplices  of the  3D 
triangulations used in the tests with finite element method \texttt{fem1} reported in Tables~\ref{pr13fem} and \ref{pr13gmsh}.
}
\label{pr13meshq}
\end{table}

\begin{table}[htbp!]\scriptsize %
\centering
\renewcommand{\arraystretch}{1.25}
\begin{tabular}{|c|c|c|c|c||c|c|}
\hline
\multirow{2}{.7cm}{$\# {\Xi _{\rm int}}$} &\multicolumn{4}{c||}{Polyharmonic $r^5$ with quadratic polynomial} & \multicolumn{2}{c|}{Pivoted QR}\\
\cline{2-7}
& {\scriptsize oct} & {\scriptsize oct-dist} & {\scriptsize 20near} & {\scriptsize pQR4sel} & {\scriptsize pQR3} & {\scriptsize pQR4}\\
\hline
815	&	5.6e{\tt-}02	&	5.1e{\tt-}02	&	1.9e{\tt-}01	&	5.0e{\tt-}02	&	4.6e{\tt-}02	&	8.4e{\tt-}01\\
1947	&	4.2e{\tt-}02	&	4.0e{\tt-}02	&	7.0e{\tt-}02	&	3.3e{\tt-}02	&	2.9e{\tt-}02	&	3.2e{\tt-}01\\
3858	&	2.6e{\tt-}02	&	2.2e{\tt-}02	&	5.0e{\tt-}02	&	2.4e{\tt-}02	&	1.8e{\tt-}02	&	6.8e{\tt-}02\\
8239	&	2.6e{\tt-}02	&	1.7e{\tt-}02	&	3.0e{\tt-}02	&	1.8e{\tt-}02	&	1.5e{\tt-}02	&	2.0e{\tt-}01\\
16386	&	1.9e{\tt-}02	&	1.4e{\tt-}02	&	2.1e{\tt-}02	&	1.4e{\tt-}02	&	1.1e{\tt-}02	&	1.3e{\tt-}01\\
35635	&	1.3e{\tt-}02	&	1.0e{\tt-}02	&	1.4e{\tt-}02	&	9.4e{\tt-}03	&	7.3e{\tt-}03	&	4.1e{\tt-}02\\
71294	&	9.2e{\tt-}03	&	7.2e{\tt-}03	&	9.9e{\tt-}03	&	6.1e{\tt-}03	&	6.0e{\tt-}03	&	6.9e{\tt-}01\\
147878	&	6.1e{\tt-}03	&	5.1e{\tt-}03	&	5.8e{\tt-}03	&	5.2e{\tt-}03	&	5.4e{\tt-}03	&	2.3e{\tt-}02\\
\hline
density	&	7.7	&	8.2	&	8.4	&	8.7	&	6.8	&	8.7\\
\hline
\end{tabular}
\caption{Test Problem~\ref{pro13}: RRMS error  $E_{\rm ref}$ for uniform interior grids and 
boundary nodes from optimized triangulations. The 7-node grid stencil \eqref{7star} is used whenever possible.
Parameters of \texttt{oct-dist}: $s=1$, $n=3$ and  $k=18$. We start with $h=0.9\,\texttt{Hmax}$ for 
$\texttt{Hmax}=H_0/2$ (that is, $i=3$) since larger $h$ led to sets $\Xi$ with too few nodes.
}
\label{pr13gridfem}
\end{table}

\begin{table}[htbp!]\scriptsize%
\centering
\renewcommand{\arraystretch}{1.25}
\begin{tabular}{|c||c|c|c||c|c|c|c||c|c|}
\hline
\multirow{2}{.7cm}{$\# {\Xi _{\rm int}}$} &\multicolumn{3}{c||}{Polyharmonic} & \multicolumn{2}{c|}{Pivoted QR}\\
\cline{2-6}
&	 {\scriptsize oct-dist}  & {\scriptsize 20near} & {\scriptsize pQR4sel}& {\scriptsize pQR3} & {\scriptsize pQR4}\\
\hline
 74	&	1.8e+01	&	8.1e-01	&	3.3e-01	&	1.9e+00	&	9.8e-01\\
174	&	1.1e+00	&	1.5e+00	&	7.5e-01	&	8.7e-01	&	7.5e+00\\
384	&	2.3e-01	&	3.8e-01	&	2.6e-01	&	9.0e-01	&	1.7e+00\\
850	&	2.0e-01	&	2.1e+00	&	7.9e-01	&	9.1e-01	&	5.3e+00\\
1828	&	1.4e-01	&	2.9e-01	&	5.8e-01	&	8.6e-01	&	3.5e+00\\
3926	&	1.0e-01	&	1.5e-01	&	3.8e+00	&	1.1e-01	&	4.4e-01\\
8238	&	5.9e-01	&	5.6e-02	&	9.1e-01	&	1.5e-01	&	2.2e+00\\
17053	&	1.6e-02	&	3.8e-02	&	1.6e-02	&	5.3e-02	&	2.8e-01\\
35193	&	1.1e-02	&	2.2e-02	&	2.0e-02	&	5.4e-02	&	2.0e-01\\
72136	&	1.4e-02	&	1.8e-02	&	2.4e-02	&	5.4e-02	&	1.0e-01\\
146898	&	5.6e-03	&	5.6e-02	&	2.9e-02	&	7.0e-02	&	5.0e-01\\
\hline
density	&	16.1	&	18.4	&	18.6	&	9.4	&	18.6\\
\hline
\end{tabular}
\caption{
Test Problem~\ref{pro13}: RRMS error $E_{\rm ref}$  for 
Halton interior nodes and projected boundary nodes.
Parameters of \texttt{oct-dist}: $s=1$, $n=3$ and  $k=17$. We skip node sets with $0\le i\le2$. %
}
\label{pr13project}
\end{table}

\subsection*{Summary observations}
\begin{enumerate}

\item Low order meshless finite difference methods with appropriate selection of the influence sets compete well 
 on presented test problems with the finite element method based on piecewise linear shape functions.

\item Existence of a shape regular triangulation of the nodes in $\Xi$ does not seem to be of importance for the mFD methods. 
This was clearly observed for Test Problem~\ref{pro11}. In other test problems the performance of mFD on 
unoptimized triangulations was significantly worse, which may be attributed to other factors such as
specific features of the triangulation algorithms applied to the STL geometry, see a discussion for Test Problem~\ref{pro19}.

\item A good discretization of the boundary was important for a good performance of mFD on the non-convex domains defined
by STL models. The best node generation methods in our experiments used boundary nodes given by the vertices of an optimized 
triangulation. For interior nodes a rather cheap but well performing alternative to an optimized 3D triangulation was 
to simply take the nodes of a Cartesian grid, which 
also has the advantage of particularly low density of the system matrix, and provides for
the possibility of using classical finite difference 7-point stencil for most interior nodes.

\item Nevertheless, all node generation methods, even the most careless, may lead to good results 
of mFD if the influence sets are well selected.

\item Stability constant $\sigma$ investigated in detail for Test Problem~\ref{pro19} seems to be a good predictor 
of the performance of the mFD methods as long as the influence sets admit polynomially exact weights \eqref{exap}. 
For node discretizations $\Xi$ generated with little care even the most reliable methods  \texttt{oct-dist} and 
\texttt{pQR4sel} are not free from occasional spikes of the stability constant that are difficult to predict, but at least 
they can be detected if $\sigma$ is computed.

\item As a cheaper alternative to the direct solvers, the solution can also be computed by the iterative method BiCGSTAB 
with ILU(0) as preconditioner whenever the mFD method performs well, as seen in the experiments for Test Problem~\ref{pro19}.

\item We cannot recommend selection method \texttt{tet} as it often fails to produce sets of influence admitting 
polynomially exact formulas. 

\item The method \texttt{oct} is significantly outperformed by \texttt{oct-dist}, whose moderate
additional cost is therefore fully justified.

\item Methods relying solely on nearest neighbors are not among the best, in particular \texttt{20near} has the highest
error on the optimized triangulation in all STL examples, and often fails on less regular nodes such as vertices of 
unoptimized triangulations. Further increasing the number of neighbors does not help, as has
been demonstrated for Test Problem~\ref{pro19}.

\item  Our new \texttt{oct-dist} method is the best all-round performer among considered selection methods, even if it is often slightly 
outperformed by \texttt{pQR3} on more regular node sets or grids. Another method deserving attention is \texttt{pQR4sel}, 
which is rather robust on all types of nodes, albeit not as much  as \texttt{oct-dist}, and typically leads to higher errors.

\item One of the features of \texttt{oct-dist} is that its parameters may be adjusted to particular types of problems and
node distributions. In particular, we used $\delta =0.9$ for Test Problems~\ref{pro11} and \ref{pro19},  
and  $\delta =0.7$ for Test Problems~\ref{pro5} and \ref{pro13} with their more complicated shapes. Parameters
$s,n,k$ were adjusted according to the type of nodes. For the nodes of unoptimized triangulations of STL models
we set $s = 3$, $n=6$, $k=17$, and in all other cases $s=1$, $n=3$, with $k$ chosen as follows. We used
 $k=13$ for optimized triangulations except of Test Problem~\ref{pro11} with smooth solution on the ball, 
 where $k=17$ was a better choice. For Cartesian grids we take $k=18$, and for Halton nodes $k=17$.

\end{enumerate}

\section{Conclusion}\label{concl}

We have introduced a new selection algorithm \texttt{oct-dist} for the Laplace operator in 3D that works reliably with 
low order RBF-FD %
on typical STL models. In particular, cheap unoptimized 3D triangulations as well as Cartesian grids or 
Halton points may be employed for 
discretization. More regular nodes obtained by optimized triangulations are beneficial but more costly to generate.
For optimized triangulations and on grids the pQR methods may be preferable.

For higher order methods the idea of the method \texttt{pQR4sel}, a combination of pQR selection with a lower order
polyharmonic RBF-FD weights may be useful. Although this method was inferior to \texttt{oct-dist} in our experiments,
it is easy to generalize to higher order, in contrast to \texttt{oct-dist}.

We have seen that any stencil selection method may fail to produce a stable system matrix in certain situations, 
hence we need a better understanding what conditions the sets of influence must satisfy in order
to generate a good system matrix. We hope that practical conditions of this type will be found in the course of the research 
on error bounds, see \cite{D19arxiv2,TLH21} for recent results that go beyond standard M-matrix arguments, 
and help to improve stencil selection algorithms. 

Our experiments for Test Problem~\ref{pro19} include estimation of the stability constant $\sigma$ defined as the 
infinity norm of the inverse of the system matrix, and solution of the system by an iterative method. They show that the size
of  $\sigma$ is a good indicator of the performance of the mFD methods, and of the convergence of the iterative solver.
Since $\sigma$ can be efficiently estimated for sparse matrices, it may be taken into account by practical algorithms for
accessing the reliability of the results.

Further experiments and possibly improvements of the selection algorithms may be needed for
practical applications, involving in particular more challenging PDEs, time-dependent settings, and more complicated CAD models, for which 
manual interventions in the mesh generation process  may be required for good results with the finite element or 
finite volume methods.

Successful experiments with the Cartesian grid in this paper and in \cite{DavySaf21} indicate that it is beneficial to use it as much as possible,
and a promising way of node generation may include combinations of Cartesian grids of various density with local improvement
of regularity of nodes in places where they meet and near boundaries and interfaces. 
Our experiments suggest that the existence of a mesh with high quality elements is not important
for a good performance of mFD, which hopefully makes node generation for mFD a significantly simpler task than mesh 
generation for mesh-based methods.

\bibliography{meshless}
\bibliographystyle{abbrv}
\end{document}